\documentclass{cas-sc}

\usepackage{ulem}
\usepackage{graphicx}
\usepackage{amsmath,amssymb}
\usepackage{xcolor}
\usepackage{bm}
\usepackage{placeins}
\usepackage{url}
\usepackage{subfigure}
\usepackage{siunitx}
\usepackage{booktabs}
\usepackage{multirow}
\usepackage{scrlayer-scrpage}

\newcommand{\pvct}[1]{\bm{#1}}
\newcommand{\vct}[1]{\bm{#1}}
\newcommand{\mtx}[1]{\pmb{\mathsf{#1}}}
\newcommand{\pxx}{\pvct{x}}
\newcommand{\pyy}{\pvct{y}}
\newcommand{\puu}{\pvct{u}}

\newtheorem{remark}{Remark}
\newtheorem{definition}{Definition}

\usepackage[capitalise]{cleveref}

\begin{document}

\title [mode = title]{A Hierarchical direct solver for Lippmann--Schwinger equations
based on boundary traces}  



%

\author[1]{Yunhui Cai}

\cormark[1]


\ead{caiyh@utexas.edu}
\shortauthors{Cai~et~al.}
\shorttitle{A~Hierarchical~direct~solver~for~Lippmann–Schwinger~equations}

\credit{Methodology, Validation, Formal analysis, Investigation, Software, Visualization, Writing -- original draft}

\affiliation[1]{organization={Department of Mathematics, The University of Texas at Austin},country={United States of America}}

\author[2]{Simon Dirckx}

\credit{Methodology, Validation, Formal analysis, Investigation, Software, Visualization, Writing -- original draft}

\affiliation[2]{organization={Oden Institute, The University of Texas at Austin},
            country={United States of America}}

\author[3]{Joar Bagge}

\credit{Methodology, Validation, Software, Writing -- original draft}

\affiliation[3]{organization={Department of Mathematics, KTH Royal Institute of Technology},
            city={SE-100\,44 Stockholm},
            country={Sweden}}

\author[1,2]{Per-Gunnar Martinsson}

\credit{Conceptualization, Methodology, Supervision, Writing -- review \& editing}

\cortext[1]{Corresponding author}



\begin{abstract}
A high-order hierarchical direct solver for variable-media scattering in two dimensions is presented.
Standard discretizations for these problems produce dense linear systems that are commonly solved iteratively, requiring a new global solve for each incident field.
The method proposed here instead constructs a reusable scattering operator for the inhomogeneous medium, based on a hierarchical subdivision of the computational domain. First, high-accuracy leaf scattering operators are built from discretized local integral equations. Up the hierarchy, neighboring boxes are merged to build the scattering operators for successively larger boxes.
A key feature of the method is that each scattering operator is represented as a map of incoming traces to outgoing traces. Consequently, volume unknowns appear only in the leaf computations; all higher levels of the hierarchy operate entirely on boundary traces.
The hierarchical solver is presented first at the continuum level and discretized only afterward using high-order local volume quadrature, polynomial trace representations, and boundary integral techniques. For a discretization involving $N$ degrees of freedom, the method requires $O(N^{3/2})$ work and $O(N\log N)$ storage to build all scattering operators. Once these are available, any incoming field can be processed at cost $O(N)$ if only the outgoing field is required, and cost $O(N\log N)$ if the full field is required. Numerical experiments demonstrate high-order accuracy on a range of challenging scattering problems including cavities, back-scattering, limited regularity, and high wavenumber, that pose significant challenges to iterative methods.
\end{abstract}


\maketitle

\section{Introduction}

\begin{figure}
\centering
\includegraphics[width=.5\linewidth]{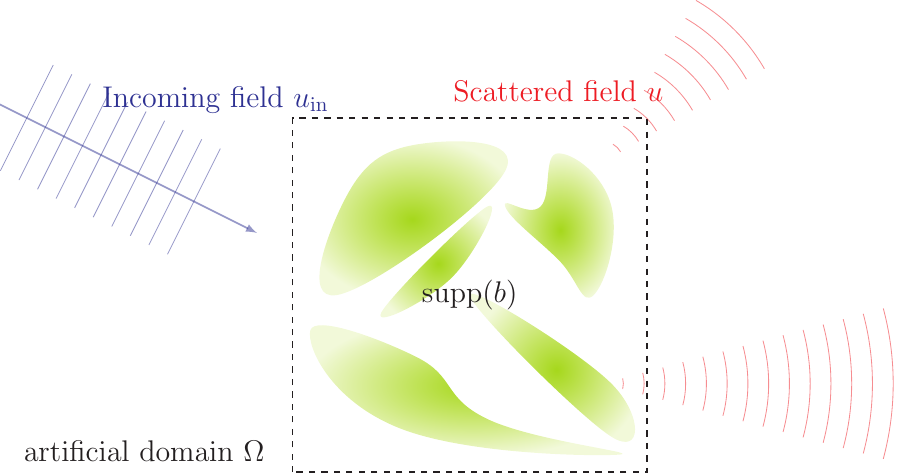}
\caption{Geometry of the scattering problem. Waves propagate in a medium with constant
wave speed $v_{\rm free}$ everywhere except in the shaded areas, where the wave speed
is given by $v(\pxx) = v_{\rm free}/\sqrt{1-b(\pxx)}$, with $b = b(\pxx)$ a given, smooth,
compactly supported ``scattering potential.'' An incident field $u_{\rm in}$ hits the
scattering potential and induces the scattered field $u$. The dashed line marks the
artificial domain $\Omega$, which encloses the support of the scattering potential.}
\label{fig:scattering_geometry}
\end{figure}

This paper describes a method for numerically solving a variable coefficient Helmholtz equation of the form
\begin{equation}
\label{eq:basic}
\left\{\begin{aligned}
-\Delta u(\pxx) - \kappa^{2}\bigl(1 - b(\pxx)\bigr)\,u(\pxx) =&\ -\kappa^{2}b(\pxx)u_{\rm in}(\pxx),
\qquad &&\mbox{for}\ \pxx \in \mathbb{R}^{2},\\
\frac{\partial u(\pxx)}{\partial r} - i\kappa u(\pxx) =&\ o(r^{-1/2}),
\qquad &&\mbox{as}\ r:=|\pxx|\rightarrow\infty.
\end{aligned}\right.
\end{equation}
Here $\kappa \in \mathbb{R}$ is the wave-number, $u_{\rm in}$ is a prescribed incoming field satisfying the free-space Helmholtz equation in a neighborhood of the scatterer, and $b$ is a compactly supported scattering potential. We let $\Omega$ denote a rectangular domain that encloses the support of $b$, and further assume for simplicity that $u_{\rm in}$ satisfies the free-space Helmholtz equation throughout $\Omega$. When $b(\pxx)<1$, an acoustic interpretation is obtained by relating $b$ to the local wave speed $v(\pxx)$ and the free-space value $v_{\rm free}$ through
$$
b(\pxx) = 1 - \left(\frac{v_{\rm free}}{v(\pxx)}\right)^{2}.
$$
The scattered field $u$ is required to satisfy the Sommerfeld radiation condition. Figure~\ref{fig:scattering_geometry} illustrates the geometry.

A particularly attractive formulation of (\ref{eq:basic}) is obtained by introducing the free-space fundamental solution
$$
G_{\kappa}(\pxx,\pyy) = \frac{i}{4}\,H_{0}^{(1)}(\kappa|\pxx-\pyy|),
$$
where $H_{0}^{(1)}$ is the Hankel function of the first kind of order zero. This is the
outgoing free-space fundamental solution for $-\Delta-\kappa^{2}I$. We then write
\begin{equation}
\label{eq:Sansatz}
u(\pxx) = \int_{\Omega} G_{\kappa}(\pxx,\pyy)\,q(\pyy)\,d\pyy,
\qquad \pxx\in\mathbb{R}^{2}.
\end{equation}
The source density $q$ then satisfies the Lippmann--Schwinger equation
\begin{equation}
\label{eq:LS2}
q(\pxx) + \kappa^{2}b(\pxx)\int_{\Omega}G_{\kappa}(\pxx,\pyy)\,q(\pyy)\,d\pyy
= -\kappa^{2}b(\pxx)\,u_{\rm in}(\pxx),
\qquad \pxx \in \Omega.
\end{equation}
The formulation (\ref{eq:LS2}) has several well-known advantages: the computational domain is finite,
the radiation condition is enforced automatically, and the equation has the
structure of a Fredholm equation of the second kind \cite{2013_colton_kress}. High-order
discretizations are available for smooth media
\cite{2004_aguilar_chen_lippmann,2008_duan_rokhlin,2016_ambikasaran_lippmann},
with related high-order volume-integral discretizations also developed for more general
acoustic and electromagnetic media
\cite{2016_anand_pandey_lippmann,2018_chen_cho_cai_vie}. 

The main computational drawback is that a discretization of (\ref{eq:LS2}) produces a
dense linear system involving all degrees of freedom in the scattering region. A standard approach to overcome this problem is to deploy an iterative solver, where each application of the operator is accelerated using a fast method such as the FFT or the Fast Multipole Method \cite{2004_aguilar_chen_lippmann,2008_duan_rokhlin}. Considerable work has been devoted to making such iterations effective at high frequency and for strongly scattering media, including sparsifying and sweeping preconditioners
\cite{2015_ying_lippmann_preconditioner,2016_zepeda_zhao_lippmann,2018_liu_ying_sparsify_sweep,2022_martinsson_lippmannschwinger}.
For many problems these techniques are highly effective and are the natural methods of
choice. There remain, however, important regimes in which the global iteration can be
difficult, and applications involving many incoming fields require repeated solution of
essentially the same scattering problem. These considerations motivate the construction
of a reusable direct solver.

The method developed in this paper is based on a simple observation. Consider any square
or rectangular subdomain $\Omega_{\tau}$. The field incident on $\Omega_{\tau}$ from
sources outside the box satisfies the homogeneous Helmholtz equation inside the box,
while the field generated by sources inside $\Omega_{\tau}$ satisfies the homogeneous
Helmholtz equation outside it. We represent both fields by their traces on the boundary
$\Gamma_{\tau}=\partial\Omega_{\tau}$. The response of the entire volume $\Omega_{\tau}$
can then be summarized by a boundary-to-boundary \textit{scattering operator}
$$
F^{\tau}:\quad \{\mbox{incoming trace on }\Gamma_{\tau}\}
\ \longmapsto\
\{\mbox{outgoing trace on }\Gamma_{\tau}\}.
$$
Once $F^{\tau}$ has been constructed, the volume degrees of freedom inside
$\Omega_{\tau}$ no longer need to participate in interactions with the rest of the problem.

The method is hierarchical. The computational domain is first partitioned into a collection
of small leaf boxes, which are organized into a binary tree by recursively pairing neighboring
boxes. Each parent box is the union of its two children, and the root of the tree is the full
computational domain $\Omega$. On each leaf, a scattering operator is constructed by solving local
Lippmann--Schwinger problems in the volume. The scattering operators are then merged
pairwise in an upward pass through the tree: given the scattering operators of two children,
a boundary-only calculation produces the scattering operator of their parent. Repeating this
process yields a scattering operator for the entire inhomogeneous medium.

A key feature of the method is that volume unknowns appear only in the leaf computations.
At all higher levels, the hierarchy operates entirely on boundary traces. This dimensional
reduction is particularly attractive at high accuracy, since inter-box interactions are
represented directly by high-order boundary fields rather than by increasingly rich subsets
of a volume discretization.

The method naturally separates into three computational stages. First, a set of
geometry- and frequency-dependent Helmholtz continuation and evaluation operators is
constructed. 
These operators depend only on the box sizes and on the wave-number $\kappa$, not on the scattering potential, and can therefore be precomputed and reused for different media.
Second, once the medium $b$ is specified, local scattering operators are built and
merged in an upward pass. Third, after this build stage, each new incoming field can be
processed rapidly through the stored hierarchy. 

For a discretization with $N$ degrees of freedom, the method requires $O(N^{3/2})$ work to construct the hierarchical scattering representation. Retaining the operators needed for subsequent recovery of the full volume solution requires $O(N\log N)$ storage, while a variant that constructs only the root scattering operator requires only $O(N)$ storage. Once the hierarchy has been built, new incident fields can be processed without any further volume solves. (Let us note that the matrices representing the scattering operators $F^{\tau}$ appear to have internal rank structure that could in principle be exploited to further reduce asymptotic complexities, see Remark \ref{rem:rank_structure}.)

In this work, the reduction to boundary traces and the hierarchical structure of scattering
operators are derived at the continuum level, before discretization. To realize the resulting
algorithm numerically, we use high-order discretizations of the local
Lippmann--Schwinger equation and high-order polynomial representations of the boundary
traces (with $p=10$ or $p=20$ being representative). The continuation and evaluation
operators are constructed using standard high-order boundary integral equation techniques.

The present approach is related to several existing lines of work. Fast iterative methods for
the Lippmann--Schwinger equation form the main computational baseline
\cite{2004_aguilar_chen_lippmann,2008_duan_rokhlin,2015_ying_lippmann_preconditioner,
2016_zepeda_zhao_lippmann,2018_liu_ying_sparsify_sweep}. Hierarchical construction of
scattering operators for volume integral equations has important precedents, notably the
work of Chen \cite{2002_chen_direct_lippman_schwinger}. T-matrix and multiple-scattering
methods likewise characterize individual scatterers by maps from incident to scattered fields
\cite{2012_ganesh_hawkins_hiptmair_tmatrix,2012_ganesh_multiple_scattering,
2019_ganesh_BIE_multiple_scattering}. There is also a substantial literature on direct
solvers that first discretize the volume integral equation and then exploit rank structure in
the resulting dense matrix
\cite{2016_ambikasaran_lippmann,corona2015n}.
A particularly close predecessor is the direct solver of \cite{2022_martinsson_lippmannschwinger}, which likewise uses a hierarchical construction based on scattering operators for the Lippmann-Schwinger equation. In \cite{2022_martinsson_lippmannschwinger}, the volume integral equation is first discretized, and rank structure in the resulting dense matrix is then exploited to construct and merge discrete scattering operators. The present construction reverses the order of these operations: the hierarchical reduction is formulated at the continuum level, so that all inter-box interactions above the leaf level are represented directly by boundary traces, with discretization introduced only afterward. This distinction becomes particularly important in the high-accuracy regime, where the interaction ranks in the discrete formulation of \cite{2022_martinsson_lippmannschwinger} can become large.

The organization of the present method is reminiscent of hierarchical Poincar\'e--Steklov (HPS) solvers for differential equations \cite{2012_spectralcomposite,2013_martinsson_DtN_linearcomplexity}:
local problems are condensed into boundary response operators and these are merged recursively. The distinction is that the local variable-coefficient scattering problems here are solved through the Lippmann--Schwinger equation rather than by directly discretizing
the differential operator.

The contribution of this paper is a high-order numerical realization of this field-to-field hierarchical scattering construction. We give the exact continuum merge identities, describe a practical discretization based on high-order local volume solves and boundary traces, and present numerical experiments on a range of variable-media scattering problems. The present implementation deliberately does not exploit the substantial numerical rank structure that
appears in many of the merge operators; this provides an immediate avenue for future acceleration.

The remainder of the paper is organized as follows. Section~\ref{sec:local_scattering}
defines the local incoming and outgoing fields and the field-to-field scattering operator.
Section~\ref{sec:merge} derives the merge operation, first for separated scatterers and then
for adjacent boxes. Section~\ref{sec:thealgorithm} describes the complete hierarchical
continuum algorithm, including the downward recovery of the interior solution.
Section~\ref{sec:disc} describes the discretization and implementation.
Numerical experiments are presented in Section~\ref{sec:num}, and conclusions are given in
Section~\ref{sec:conc}.

\section{Local scattering operators}
\label{sec:local_scattering}

\begin{figure}
\centering
\includegraphics[width=.5\linewidth]{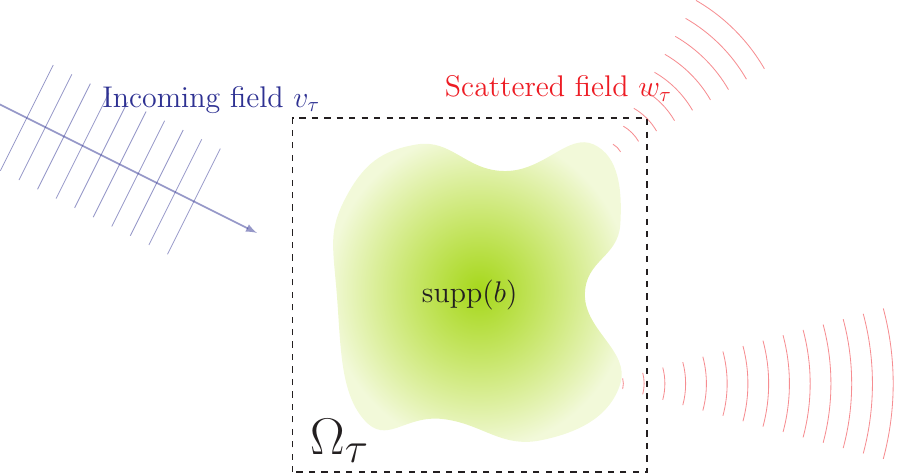}
\caption{A box $\Omega_{\tau}$ contains the support of a scattering potential $b$.
An incoming field $v_{\tau}$ induces the scattered field $w_{\tau}$. The fields are represented by their traces $\phi_{\tau}$ and $\psi_{\tau}$ on the artificial boundary
$\Gamma_{\tau}=\partial\Omega_{\tau}$, so the local scattering operator $F^{\tau}$ maps $\psi_{\tau}=F^{\tau}\phi_{\tau}$.}
\label{fig:singlebox}
\end{figure}

In this section, we give a precise mathematical definition of the  local scattering operators that will serve as the basic building blocks of our hierarchical solver. 
Throughout this section, $\tau$ denotes a label associated with a square or rectangular subdomain $\Omega_\tau\subset\Omega$, with boundary $\Gamma_\tau=\partial\Omega_\tau$, cf.~Figure  \ref{fig:singlebox}. These subdomains will later form part of the hierarchical tree.

For a single square or rectangular subdomain, we first identify the incoming and outgoing fields, then represent these fields by their boundary traces, and finally show how the resulting boundary-to-boundary scattering operator can be constructed from a local Lippmann--Schwinger solve.

\subsection{Incoming and outgoing fields}

For a box $\Omega_{\tau}$, we define \textit{incoming} and
\textit{outgoing} fields, following standard practice. Let $q$ denote the solution to
(\ref{eq:LS2}).
The \textit{outgoing field} generated by the sources in $\Omega_{\tau}$ is
\begin{equation}
\label{eq:def_outgoing}
w_{\tau}(\pxx) = \int_{\Omega_{\tau}}G_{\kappa}(\pxx,\pyy)\,q(\pyy)\,d\pyy,
\qquad \pxx \in \mathbb{R}^{2}\setminus\overline{\Omega}_{\tau}.
\end{equation}
In other words, $w_{\tau}$ is the radiating solution to the $\kappa$-Helmholtz equation generated by the sources inside $\Omega_{\tau}$.

The \textit{incoming field} for the box $\tau$ is the sum of the externally applied field
$u_{\rm in}$ and the field generated by all sources outside $\Omega_{\tau}$. Thus,
\begin{equation}
\label{eq:def_incoming}
v_{\tau}(\pxx) = u_{\rm in}(\pxx) +
\int_{\Omega\setminus\Omega_{\tau}}G_{\kappa}(\pxx,\pyy)\,q(\pyy)\,d\pyy,
\qquad \pxx \in \Omega_{\tau}.
\end{equation}
Observe that $v_{\tau}$ satisfies the homogeneous Helmholtz equation throughout
$\Omega_{\tau}$.

Once the incoming field $v_{\tau}$ is known, the Lippmann--Schwinger equation can be
solved locally on $\Omega_{\tau}$. Indeed, restricting (\ref{eq:LS2}) to
$\Omega_{\tau}$ and splitting the volume potential into an interior and an exterior
part gives
\begin{equation}
\label{eq:LS2loc}
q(\pxx) + \kappa^{2}b(\pxx)
\int_{\Omega_{\tau}}G_{\kappa}(\pxx,\pyy)\,q(\pyy)\,d\pyy
=
-\kappa^{2}b(\pxx)v_{\tau}(\pxx),
\qquad \pxx \in \Omega_{\tau}.
\end{equation}
Solving (\ref{eq:LS2loc}) specifies the restriction of $q$ to $\Omega_{\tau}$,
and hence the outgoing field $w_{\tau}$ via (\ref{eq:def_outgoing}).

\subsection{Trace representations}
\label{sec:scatter_trace}

When a scattering operator is to be computed in practice, one must specify how the
incoming and outgoing fields are represented. Standard practice in the physics literature
is to use truncated expansions in cylindrical harmonics (and, in three dimensions,
spherical harmonics). We will instead
represent both fields through their \textit{boundary traces}, by which we simply mean
their restrictions to the boundary. For a box $\tau$, define
\begin{equation}
\label{eq:def_traces}
\phi_{\tau}=v_{\tau}|_{\Gamma_{\tau}},
\qquad
\psi_{\tau}=w_{\tau}|_{\Gamma_{\tau}}.
\end{equation}
Thus, $\phi_{\tau}$ is the incoming trace and $\psi_{\tau}$ is the outgoing trace. 
(The incoming and outgoing fields considered here are continuous, so 
the boundary traces are well defined.)

The outgoing trace $\psi_{\tau}$ always uniquely determines $w_{\tau}$: it supplies
Dirichlet data for the exterior Helmholtz problem, which is uniquely solvable when coupled
with the radiation condition at infinity. The incoming trace $\phi_{\tau}$ supplies
Dirichlet data for an interior Helmholtz problem on $\Omega_{\tau}$. For Helmholtz,
this problem can in principle be ill-posed. For the present development, make the assumption:
\begin{center}
\fbox{For every box $\Omega_{\tau}$ in the hierarchy, the local Dirichlet problem is well posed.}
\end{center}
Under this assumption, $\phi_{\tau}$ uniquely determines the incoming field
$v_{\tau}$ throughout $\Omega_{\tau}$.\footnote{It may be possible to relax this
assumption by switching to an impedance trace formulation. This is a subject of
future work. Observe, however, that the potentially dangerous resonances are for
the free-space Helmholtz equation in rectangular boxes. They occur when $\kappa^{2}$
coincides with a Dirichlet eigenvalue of $-\Delta$, and these eigenvalues are known
analytically. Since we can choose the box sizes, it is in principle easy to avoid exact
resonances. Ill-conditioning near a resonance is harder to avoid, however.}

We now define the boundary-to-boundary scattering operator $F^{\tau}$ by
\begin{equation}
\label{eq:def_scattering_operator}
F^{\tau}\,\colon\,\phi_{\tau}\mapsto\psi_{\tau}.
\end{equation}
The scattering operator encodes the complete response of the material in
$\Omega_{\tau}$ to any incoming field that is homogeneous Helmholtz in that box.

\subsection{Constructing a local scattering operator}
\label{sec:leaf_comp_continuum}

The scattering operator for a leaf box can be written as a composition of four continuum
operators. Let
\begin{itemize}
\item $S_{\tau}$ be the solution operator for the local Dirichlet Helmholtz problem on
$\Omega_{\tau}$, so that
$$
v_{\tau}=S_{\tau}\phi_{\tau}.
$$
\item $B_{\tau}$ be multiplication by the function $-\kappa^{2}b(\pxx)$.
\item $L_{\tau}$ be the locally restricted Lippmann--Schwinger operator
$$
[L_{\tau}q](\pxx)=q(\pxx)+\kappa^{2}b(\pxx)
\int_{\Omega_{\tau}}G_{\kappa}(\pxx,\pyy)\,q(\pyy)\,d\pyy,
\qquad \pxx\in\Omega_{\tau}.
$$
\item $M_{\tau}$ evaluate the outgoing field on the boundary:
$$
[M_{\tau}q](\pxx)=
\int_{\Omega_{\tau}}G_{\kappa}(\pxx,\pyy)\,q(\pyy)\,d\pyy,
\qquad \pxx\in\Gamma_{\tau}.
$$
\end{itemize}
Define the local source operator
\begin{equation}
\label{eq:def_source_operator}
Q^{\tau}=L_{\tau}^{-1}B_{\tau}S_{\tau}.
\end{equation}
Then $q|_{\Omega_{\tau}}=Q^{\tau}\phi_{\tau}$, and the leaf scattering operator is
\begin{equation}
\label{eq:leaf_scattering_operator}
F^{\tau}=M_{\tau}Q^{\tau}
=M_{\tau}L_{\tau}^{-1}B_{\tau}S_{\tau}.
\end{equation}
Formula (\ref{eq:leaf_scattering_operator}) will be used only on leaf boxes. For larger
boxes, the scattering operator will instead be constructed by merging the operators of
two smaller boxes.

\begin{remark}[Equivalent-source formulation]
\label{rem:equivalent_source}
The field-to-field scattering operator used throughout the main text is not the only
possible boundary representation. An alternative is to continue representing the
incoming field by its trace $\phi_{\tau}$, but to represent the outgoing field by an
equivalent single-layer source density supported on $\Gamma_{\tau}$. This leads to a
closely related hierarchical merge in which interactions between children are evaluated
directly through single-layer potentials. We use the field-to-field formulation in the
main text because it proved simpler to implement and avoids the reconstruction of
equivalent source densities. The alternative formulation is recorded in Appendix~\ref{app:equivalent_source}.
\end{remark}

\section{Merging scattering operators}
\label{sec:merge}

Suppose that a subdomain $\Omega_{\tau}$ is formed as a union $\Omega_{\tau} = \Omega_{\alpha} \cup \Omega_{\beta}$ of two smaller subdomains $\Omega_{\alpha}$ and $\Omega_{\beta}$,
and that the scattering operators $F^{\alpha}$ and $F^{\beta}$ for the two ``children'' have already been constructed. We will show how to form the scattering operator $F^{\tau}$ for the ``parent'' using boundary-to-boundary operations only. Once the scattering operator of a box has been built, we never need to revisit its interior during the upward pass.

To introduce the merge in the cleanest possible setting, we first consider two children
that are separated in space. We then turn to the case of two adjacent boxes, which is the
case required by the hierarchical solver. The second derivation involves additional
bookkeeping, but the underlying algebra is the same.

\subsection{Merging two separated boxes}
\label{sec:mergetwosep}

\begin{figure}
\begin{center}
\begin{minipage}{.48\linewidth}
\includegraphics[width=\linewidth]{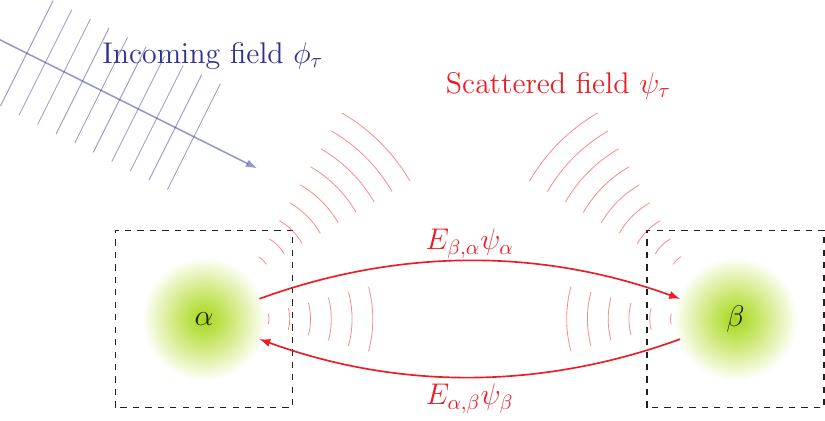}    
\end{minipage}
\hfill
\fbox{\begin{minipage}{.4\linewidth}
\small
Coupling equation:
$$
\left[\begin{array}{r} \phi_{\alpha} \\ \phi_{\beta} \end{array}\right] =
\phi_{\tau} +
\left[\begin{array}{rr} 0 & E_{\alpha,\beta} \\ E_{\beta,\alpha} & 0 \end{array}\right]
\left[\begin{array}{r} \psi_{\alpha} \\ \psi_{\beta} \end{array}\right]
$$
Child scattering:
$$
\left[\begin{array}{r} \psi_{\alpha} \\ \psi_{\beta} \end{array}\right] =
\left[\begin{array}{rr} F^{\alpha} & 0 \\ 0 & F^{\beta} \end{array}\right]
\left[\begin{array}{r} \phi_{\alpha} \\ \phi_{\beta} \end{array}\right]
$$
\end{minipage}}
\end{center}
\caption{Merging of two scattering operators, as discussed in Section \ref{sec:mergetwosep}.
Let $\Omega_{\alpha}$ and $\Omega_{\beta}$ be two boxes, with the union $\Omega_{\tau} =
\Omega_{\alpha} \cup \Omega_{\beta}$. If the scattering operators $F^{\alpha}$ and $F^{\beta}$
have been constructed, then the joint scattering operator $F^{\tau}$ can be constructed without
the need to revisit the interiors of $\Omega_{\alpha}$ or $\Omega_{\beta}$; only operators
supported on the boundaries are involved. The parent operator maps the incoming trace
$\phi_{\tau}$ to the outgoing trace $\psi_{\tau}$.}
\label{fig:twoboxes_sep}
\end{figure}

Consider a scattering potential supported in two disjoint boxes $\Omega_{\alpha}$ and
$\Omega_{\beta}$, and suppose that $F^{\alpha}$ and $F^{\beta}$ have been built,
cf.~Figure \ref{fig:twoboxes_sep}. Set
$$
\Omega_{\tau}=\Omega_{\alpha}\cup\Omega_{\beta},
\qquad
\Gamma_{\tau}=\Gamma_{\alpha}\cup\Gamma_{\beta}.
$$
Since the boxes are separated, the boundary of the parent is simply the disjoint union of the two child boundaries. (For this introductory derivation, we temporarily allow the parent domain $\Omega_\tau$ to be disconnected; the hierarchical solver itself will use only adjacent rectangular children.)

Let $\phi_{\tau}$ denote the trace on $\Gamma_{\tau}$ of an externally applied incoming field. Since $\Gamma_{\tau}=\Gamma_{\alpha}\cup\Gamma_{\beta}$, we regard $\phi_{\tau}$ as a two-component trace ordered by $(\Gamma_{\alpha},\Gamma_{\beta})$. Note that $\phi_{\tau}\neq [\phi_{\alpha},\phi_{\beta}]^T$, because in the latter quantities we must also consider the interaction between $\alpha$ and $\beta$. For later convenience (see Section~\ref{sec:mergetwo}), we introduce
$$
\eta=
\left[\begin{array}{c}
\phi_{\alpha}\\
\phi_{\beta}
\end{array}\right]
\quad\text{and}\quad
\gamma=
\left[\begin{array}{c}
\psi_{\alpha}\\
\psi_{\beta}
\end{array}\right].
$$
Introduce the block-diagonal child scattering operator
\begin{equation}
\label{eq:F_separated}
F=
\left[\begin{array}{cc}
F^{\alpha}&0\\
0&F^{\beta}
\end{array}\right],
\end{equation}
so that
\begin{equation}
\label{eq:sep_scattering}
\gamma=F\eta.
\end{equation}

The sibling-generated incoming field is obtained by propagating each outgoing field to the
other box. Let $E_{\alpha,\beta}$ denote the operator that takes the outgoing trace on
$\Gamma_{\beta}$, reconstructs the corresponding radiating exterior Helmholtz field,
and evaluates this field on $\Gamma_{\alpha}$. Define $E_{\beta,\alpha}$ analogously,
and set
\begin{equation}
\label{eq:E_separated}
E_{\tau}=
\left[\begin{array}{cc}
0&E_{\alpha,\beta}\\
E_{\beta,\alpha}&0
\end{array}\right].
\end{equation}
The total incoming trace on each child is the sum of the externally applied field and the
field generated by its sibling. Thus
\begin{equation}
\label{eq:sep_sibling_incoming}
\eta=\phi_{\tau}+E_{\tau}\gamma.
\end{equation}
Combining (\ref{eq:sep_scattering}) and (\ref{eq:sep_sibling_incoming}), we obtain the \textit{equilibrium equation} for the local incoming fields
\begin{equation*}
\left[\begin{array}{r} \phi_{\alpha} \\ \phi_{\beta} \end{array}\right] =
\phi_{\tau} +
\left[\begin{array}{cc} 0 & E_{\alpha,\beta} \\ E_{\beta,\alpha} & 0 \end{array}\right]
\left[\begin{array}{cc} F_{\alpha}&0 \\ 0&F_{\beta} \end{array}\right]
\left[\begin{array}{r} \phi_{\alpha} \\ \phi_{\beta} \end{array}\right]
\end{equation*}
or, more compactly,
$$
\eta=\phi_{\tau}+E_{\tau}F\eta.
$$
This motivates the following definition.

\begin{definition}[Solution-transfer operator]
\label{def:solution_transfer_separated}
The solution-transfer operator for the separated pair is
\begin{equation}
\label{eq:T_separated}
T^{\tau}=(I-E_{\tau}F)^{-1}.
\end{equation}
It maps the incoming trace on the parent to the total incoming traces experienced by
the two children:
\begin{equation}
\label{eq:sep_child_incoming}
\eta=T^{\tau}\phi_{\tau}.
\end{equation}
\end{definition}

We assume that the local Dirichlet problems and the parent scattering problem are uniquely solvable, 
so that the inverse in the merge formula (\ref{eq:T_separated}) exists. 

It remains to form the outgoing trace $\psi_{\tau}$ of the composite scatterer. On
$\Gamma_{\alpha}$, the total outgoing field is the sum of the field generated by
$\alpha$ and the field generated by $\beta$; the analogous statement holds on
$\Gamma_{\beta}$. Consequently,
\begin{equation}
\label{eq:sep_parent_outgoing}
\psi_{\tau}=(I+E_{\tau})\gamma.
\end{equation}
Combining (\ref{eq:sep_scattering}), (\ref{eq:sep_child_incoming}), and
(\ref{eq:sep_parent_outgoing}), we obtain
\begin{equation}
\label{eq:sep_merge}
\boxed{
T^{\tau}=(I-E_{\tau}F)^{-1},
\qquad
F^{\tau}=(I+E_{\tau})F T^{\tau}.
}
\end{equation}
Thus the parent scattering operator and the information needed to recover the fields on
its children are both constructed using boundary-to-boundary operators only.

\begin{remark}[Multiple scattering]
\label{rem:multiple_scattering}
The solution-transfer operator (\ref{eq:T_separated}) has a natural physical interpretation
in terms of repeated scattering between the children. Observe that
\begin{equation}
\label{eq:T_neumann}
T^{\tau}
=
(I - E_{\tau}F)^{-1}
=
I+E_{\tau}F+(E_{\tau}F)^{2}+(E_{\tau}F)^{3}+\cdots,
\end{equation}
whenever the series converges.
The term $E_{\tau}F$ represents one scattering event followed by propagation to the sibling,
$(E_{\tau}F)^{2}$ represents two successive interactions, and so forth. Thus
$T^{\tau}$ incorporates all orders of multiple scattering between the children.
\end{remark}

\subsection{Merging two adjacent boxes}
\label{sec:mergetwo}
\begin{figure}
\begin{center}
\includegraphics[height=30mm]{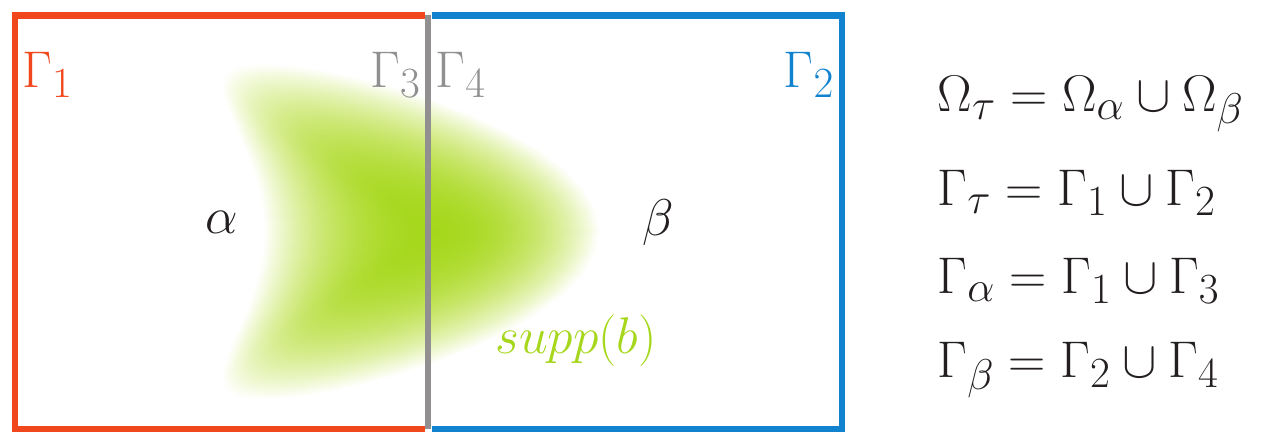}
\end{center}
\caption{Two adjacent boxes to be merged. The labels $1,2,3,4$ refer to boundary
pieces. The labels $3$ and $4$ denote two copies of the same artificial interface, so that
the two children may carry distinct trace variables and, after discretization, distinct grids.
The support of $b$ may cross this interface.}
\label{fig:twoboxes}
\end{figure}

Consider next a case where the two boxes $\alpha$ and $\beta$ directly touch,
as illustrated in Figure \ref{fig:twoboxes}.
This is the situation of immediate relevance to our hierarchical solution process, since
it consists of pairwise merging of domains to build scattering operators for larger and
larger subdomains. Let
$$
\partial\Omega_{\alpha}=\Gamma_{1}\cup\Gamma_{3},
\qquad
\partial\Omega_{\beta}=\Gamma_{2}\cup\Gamma_{4},
$$
where $\Gamma_{3}=\Gamma_{4}$ is the common interface. The parent box
$\Omega_{\tau}=\Omega_{\alpha}\cup\Omega_{\beta}$ has boundary
$$
\Gamma_{\tau}=\Gamma_{1}\cup\Gamma_{2}.
$$
We use two labels for the common interface because its traces belong to different child
boxes, and because the two copies may later be discretized differently.

The adjacent case differs from the separated case in three ways. First, the externally
applied incoming trace on the parent boundary does not include its values on the common
interface, so these values must be reconstructed. Second, the sibling-evaluation operator
contains identity blocks on the common interface; this makes computations cheaper, but
requires additional notational machinery. Third, the child outgoing traces include the
common interface, which must be eliminated when the outgoing trace of the parent is formed.

\subsubsection{Continuing the externally applied incoming field to the child boundaries}

Let $\phi_{\tau}$ denote the incoming trace on the parent boundary $\Gamma_{\tau}$.
The corresponding field satisfies the homogeneous Helmholtz equation throughout
$\Omega_{\tau}$, so its trace on the common interface can be reconstructed by solving the
interior Dirichlet problem on $\Omega_{\tau}$. We define the incoming continuation operator
$U_{\tau}$ to map $\phi_{\tau}$ to the trace of this same field on the complete boundaries
of the two children. In the boundary ordering $(1,2,3,4)$, this operator has the form
\begin{equation}
\label{eq:U_structure}
U_{\tau}=
\left[\begin{array}{cc}
I&0\\
0&I\\
U_{31}&U_{32}\\
U_{41}&U_{42}
\end{array}\right],
\end{equation}
so that $U_{\tau}\phi_{\tau}$ gives the trace of the externally applied incoming field on
the complete boundaries of the two children. Since $\Gamma_{3}=\Gamma_{4}$, the third and
fourth components of $U_{\tau}\phi_{\tau}$ agree at the continuum level; equivalently,
$U_{31}=U_{41}$ and $U_{32}=U_{42}$.

\subsubsection{Resolving the interaction between the children}

We now split the scattering operators of the children according to their exterior and
interface boundary pieces:
\begin{equation}
\label{eq:F_child_blocks}
F^{\alpha}=
\left[\begin{array}{cc}
F_{11}^{\alpha}&F_{13}^{\alpha}\\
F_{31}^{\alpha}&F_{33}^{\alpha}
\end{array}\right],
\qquad
F^{\beta}=
\left[\begin{array}{cc}
F_{22}^{\beta}&F_{24}^{\beta}\\
F_{42}^{\beta}&F_{44}^{\beta}
\end{array}\right].
\end{equation}
In the boundary ordering $(1,2,3,4)$, define
\begin{equation}
\label{eq:F_adjacent_abstract}
F=
\left[\begin{array}{cccc}
F_{11}^{\alpha}&0&F_{13}^{\alpha}&0\\
0&F_{22}^{\beta}&0&F_{24}^{\beta}\\
F_{31}^{\alpha}&0&F_{33}^{\alpha}&0\\
0&F_{42}^{\beta}&0&F_{44}^{\beta}
\end{array}\right].
\end{equation}
Thus $F$ simply applies the scattering operator of each child to the appropriate pair of
boundary traces.
With the ordering introduced above, the total incoming and outgoing traces $\eta,\gamma$ on the two children from Section~\ref{sec:mergetwosep} become
$$
\eta=
\left[\begin{array}{c}
\eta_{1}\\\eta_{2}\\\eta_{3}\\\eta_{4}
\end{array}\right]\quad \text{and} \quad
\gamma=
\left[\begin{array}{c}
\gamma_{1}\\\gamma_{2}\\\gamma_{3}\\\gamma_{4}
\end{array}\right],
$$
with 
$$
\phi_{\alpha}=\left[\begin{array}{c}\eta_{1}\\\eta_{3}\end{array}\right],
\qquad
\phi_{\beta}=\left[\begin{array}{c}\eta_{2}\\\eta_{4}\end{array}\right],
\qquad
\psi_{\alpha}=\left[\begin{array}{c}\gamma_{1}\\\gamma_{3}\end{array}\right],
\qquad
\psi_{\beta}=\left[\begin{array}{c}\gamma_{2}\\\gamma_{4}\end{array}\right].
$$
The child scattering relations are simply, cf.~(\ref{eq:sep_scattering}),
\begin{equation}
\label{eq:adj_scattering}
\gamma=F\eta.
\end{equation}

Define $E_{\tau}$ to be the operator that maps the child outgoing traces to the incoming
traces generated on the sibling. On the common interface, no exterior solve is required:
the outgoing trace from one child is directly the incoming trace generated on the other
child. In the ordering $(1,2,3,4)$,
\begin{equation}
\label{eq:E_adjacent}
E_{\tau}=
\left[\begin{array}{cccc}
0&E_{12}&0&E_{14}\\
E_{21}&0&E_{23}&0\\
0&0&0&I\\
0&0&I&0
\end{array}\right].
\end{equation}
The total incoming trace on each child is the sum of the parent incoming field, continued to
the child boundary, and the field generated by its sibling. Thus, cf.~(\ref{eq:sep_sibling_incoming}),
\begin{equation}
\label{eq:adj_sibling_incoming}
\eta=U_{\tau}\phi_{\tau}+E_{\tau}\gamma.
\end{equation}
Combining (\ref{eq:adj_scattering}) and (\ref{eq:adj_sibling_incoming}), we find
$$
\eta=U_{\tau}\phi_{\tau}+E_{\tau}F\eta,
$$
or
$$
(I-E_{\tau}F)\eta=U_{\tau}\phi_{\tau}.
$$
We define the solution-transfer operator for the parent box by
\begin{equation}
\label{eq:T_adjacent}
\boxed{
T^{\tau}=(I-E_{\tau}F)^{-1}U_{\tau}.
}
\end{equation}
Then
\begin{equation}
\label{eq:adj_child_incoming}
\eta=T^{\tau}\phi_{\tau}.
\end{equation}
The operator $T^{\tau}$ will be used twice: it enters the construction of the parent
scattering operator during the upward pass, and it transfers the incoming trace on a
parent to the total incoming traces on its children during the downward pass.

\subsubsection{Forming the parent outgoing field}

The outgoing traces in $\gamma$ are tabulated on the complete child boundaries, including
the common interface. The outgoing trace of the parent is required only on
$\Gamma_{\tau}=\Gamma_{1}\cup\Gamma_{2}$ and must include the fields generated by
both children. We therefore introduce the outgoing assembly operator $R_{\tau}$ by
\begin{equation}
\label{eq:def_Rtau}
\psi_{\tau}=R_{\tau}\gamma.
\end{equation}
In the ordering $(1,2,3,4)$,
\begin{equation}
\label{eq:R_structure}
R_{\tau}=
\left[\begin{array}{cccc}
I&E_{12}&0&E_{14}\\
E_{21}&I&E_{23}&0
\end{array}\right].
\end{equation}
Combining (\ref{eq:adj_scattering}), (\ref{eq:adj_child_incoming}), and
(\ref{eq:def_Rtau}), we obtain the parent scattering operator
\begin{equation}
\label{eq:adj_merge}
\boxed{
F^{\tau}=R_{\tau}F T^{\tau},
\qquad
T^{\tau}=(I-E_{\tau}F)^{-1}U_{\tau}.
}
\end{equation}
These are the two fundamental merge formulas.

\subsection{Summary of the merge}
\label{sec:merge_summary}

The separated and adjacent cases are governed by the same three operations:
\begin{enumerate}
\item continue the externally applied incoming field on the parent to the complete child boundaries;
\item resolve all multiple scattering between the children;
\item assemble the child outgoing fields into the outgoing field of the parent.
\end{enumerate}
The first and third operations are represented by $U_{\tau}$ and $R_{\tau}$,
respectively, while the first two combine to form the solution-transfer operator
$T^{\tau}$. If $F$ denotes the block operator formed from the two child scattering
operators, then the generic merge is
\begin{equation}
\label{eq:generic_merge}
\boxed{
T^{\tau}=(I-E_{\tau}F)^{-1}U_{\tau},
\qquad
F^{\tau}=R_{\tau}F T^{\tau}.
}
\end{equation}
For separated scatterers, $E_{\tau}$ is the operator defined in
(\ref{eq:E_separated}), $U_{\tau}=I$, and $R_{\tau}=I+E_{\tau}$. For adjacent
boxes, $U_{\tau}$ reconstructs the externally applied incoming field on the internal
interface, $E_{\tau}$ contains identity maps across that interface, and $R_{\tau}$
eliminates the interface while accumulating both child contributions on the parent boundary.
No approximation has entered the construction in this section.

\begin{remark}[Universal operators]\label{rem:universal_ops}
The operators $U_{\tau}$, $E_{\tau}$, and $R_{\tau}$ are independent of the scattering
potential $b$. They depend only on the box geometry and the free-space Helmholtz operator.
Likewise, the leaf operators $S_{\tau}$ and $M_{\tau}$ are universal. All dependence on
$b$ enters through $B_{\tau}$ and $L_{\tau}$ in the local Lippmann--Schwinger solve, and
hence through $Q^{\tau}$, $F^{\tau}$, and $T^{\tau}$. Consequently, the universal operators
can be precomputed and reused for different scattering potentials. On the uniform,
translation-invariant trees considered here, only one copy is needed for each box size and
merge orientation.
\end{remark}

\begin{remark}[Equivalent formulation and rank structure]
\label{rem:rank_structure}
The identity
$$
F(I-E_{\tau}F)^{-1}=(I-FE_{\tau})^{-1}F
$$
shows that the parent scattering operator can equivalently be written as
\begin{equation}
\label{eq:merge_outgoing_form}
F^{\tau}=R_{\tau}(I-FE_{\tau})^{-1}F U_{\tau}.
\end{equation}
This form is useful when considering possible accelerations of the merge. From
(\ref{eq:F_adjacent_abstract}) and (\ref{eq:E_adjacent}), we obtain
\begin{equation}
\label{eq:FE_full_structure}
FE_{\tau}=
\left[\begin{array}{c|c|c|c}
0&F_{11}^{\alpha}E_{12}&0&F_{11}^{\alpha}E_{14}+F_{13}^{\alpha}\\ \hline
F_{22}^{\beta}E_{21}&0&F_{22}^{\beta}E_{23}+F_{24}^{\beta}&0\\ \hline
0&F_{31}^{\alpha}E_{12}&0&F_{31}^{\alpha}E_{14}+F_{33}^{\alpha}\\ \hline
F_{42}^{\beta}E_{21}&0&F_{42}^{\beta}E_{23}+F_{44}^{\beta}&0
\end{array}\right].
\end{equation}
For boxes that contain a moderate number of wavelengths, numerical experiments indicate
that all nonzero blocks in (\ref{eq:FE_full_structure}) except the $(3,4)$ and $(4,3)$
blocks have low numerical rank. The two exceptional blocks describe interactions across
the common interface. They are not globally low rank, but admit hierarchically
semiseparable (HSS) representations. We do not pursue this structure in the present work. Equation (\ref{eq:FE_full_structure}) nevertheless points to
a particularly simple avenue for accelerating the method in future work. As the number of
wavelengths in a box grows, the numerical ranks naturally grow as well.
\end{remark}

\section{The hierarchical continuum algorithm}
\label{sec:thealgorithm}

We are now ready to describe the complete continuum algorithm for solving the
Lippmann--Schwinger equation (\ref{eq:LS2}). The method is hierarchical. It consists of
partitioning the computational domain into small leaf boxes, constructing a local scattering
operator on every leaf, and then repeatedly merging pairs of child operators in an upward
pass through a binary tree. The solution-transfer operators constructed during the upward
pass are subsequently used to recover the incoming traces on all leaves in a downward pass.

\subsection{The hierarchical partition}

\begin{figure}
\includegraphics[width=\textwidth]{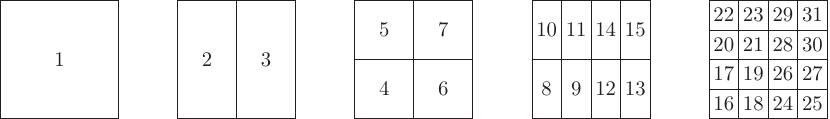}
\caption{The square domain $\Omega$ is split into $4 \times 4$ leaf boxes.
These are then gathered into a binary tree of successively larger boxes
as described in Section \ref{sec:thealgorithm}. One possible enumeration
of the boxes in the tree is shown, but note that the only restriction is
that if box $\tau$ is the parent of box $\sigma$, then $\tau < \sigma$.}
\label{fig:tree_numbering}
\end{figure}

The computational domain $\Omega$ will always be a square or rectangle. Observe that the
domain is in this context an artificial boundary box anyway. We first partition $\Omega$
into a collection of square, or possibly rectangular, leaf boxes. These should be small
enough that the Lippmann--Schwinger equation restricted to a leaf can easily be solved to
very high accuracy. We then construct a binary tree by hierarchically merging the boxes,
making sure that all boxes on the same level are roughly of the same size; see Figure
\ref{fig:tree_numbering}. The boxes are ordered so that if $\tau$ is the parent of
$\sigma$, then $\tau<\sigma$. The root of the tree, which represents the full domain
$\Omega$, has index $\tau=1$.

\subsection{Upward pass}

Let $N_{\rm boxes}$ denote the total number of boxes in the tree. The upward pass
constructs the scattering operator $F^{\tau}$ for every box and the solution-transfer
operator $T^{\tau}$ for every nonleaf box. At a leaf, we also retain the
local source operator $Q^{\tau}$ from (\ref{eq:def_source_operator}), since it will be
needed to recover the volume density after the incoming trace on the leaf has been found.

\begin{figure}
\centering
\fbox{\begin{minipage}[t][4.2cm][t]{.5\linewidth}
\begin{tabbing}
\mbox{}\hspace{6mm} \= \mbox{}\hspace{6mm} \= \mbox{}\hspace{6mm} \= \mbox{}\hspace{6mm} \= \kill
\textbf{for} $\tau = N_{\rm boxes},\,N_{\rm boxes}-1,\,N_{\rm boxes}-2,\,\dots,\,1$\\
\> \textbf{if} ($\tau$ is a leaf)\\
\> \> Build $Q^{\tau}=L_{\tau}^{-1}B_{\tau}S_{\tau}$.\\
\> \> Build $F^{\tau}=M_{\tau}Q^{\tau}$.\\
\> \textbf{else}\\
\> \> Let $\{\alpha,\beta\}$ be the children of $\tau$.\\
\> \> Form $F$ from $F^{\alpha}$ and $F^{\beta}$ as in (\ref{eq:F_adjacent_abstract}).\\
\> \> Build $T^{\tau}=(I-E_{\tau}F)^{-1}U_{\tau}$.\\
\> \> Build $F^{\tau}=R_{\tau}F T^{\tau}$.\\
\> \textbf{end if}\\
\textbf{end for}
\end{tabbing}
\end{minipage}}.
\fbox{\begin{minipage}[t][4.2cm][t]{.5\linewidth}
\begin{tabbing}
\mbox{}\hspace{6mm} \= \mbox{}\hspace{6mm} \= \mbox{}\hspace{6mm} \= \mbox{}\hspace{6mm} \= \kill
\textbf{initialize} $\phi_1 = u_{{\rm in}|\partial\Omega}$\\
\textbf{for} $\tau = 1,\ldots,N_{\rm boxes}$\\
\> \textbf{if} ($\tau$ is a leaf)\\
\> \> Compute $q_{\tau} = Q^{\tau}\phi_{\tau}$.\\
\> \> Compute $u|_{\Omega_\tau}$ as in \eqref{eq:recover_leaf_incoming}--\eqref{eq:scattered_field_in_leaf}.\\
\> \textbf{else}\\
\> \> Let $\{\alpha,\beta\}$ be the children of $\tau$.\\
\> \> Compute $\eta = T^{\tau}\phi_{\tau}$.\\
\> \> Set $\phi_{\alpha}\leftarrow\begin{bmatrix}
    \eta_1& \eta_3
\end{bmatrix}^T$ and $\phi_{\beta}\leftarrow\begin{bmatrix}
    \eta_2& \eta_4
\end{bmatrix}^T$.\\
\> \textbf{end if}\\
\textbf{end for}
\end{tabbing}
\end{minipage}}
\label{fig:upanddown}
\caption{Upwards (left) and downward (right) pass algorithms for the construction of the global scattering operator and the computation of the interior traces.}
\end{figure}

The root operator $F^{1}$ maps the incoming trace on $\partial\Omega$ to the outgoing
trace generated by all sources in $\Omega$. In particular, it determines the complete
exterior scattering solution without ever revisiting the interiors of the boxes.

\subsection{Downward pass and recovery of the solution}\label{sec:downward-pass}

Let the incoming trace on the root be
\begin{equation}
\label{eq:root_incoming_trace}
\phi_{1}=u_{\rm in}|_{\partial\Omega}.
\end{equation}
Suppose that the incoming trace $\phi_{\tau}$ has been determined for a nonleaf box
$\tau$ with children $\alpha$ and $\beta$. By the definition of the solution-transfer
operator,
\begin{equation}
\label{eq:downward_transfer}
\eta=T^{\tau}\phi_{\tau}.
\end{equation}
For the ordering $(1,2,3,4)$ used in Section~\ref{sec:mergetwo}, the incoming traces on
the two children are then
$$
\phi_{\alpha}=\left[\begin{array}{c}\eta_{1}\\\eta_{3}\end{array}\right],
\qquad
\phi_{\beta}=\left[\begin{array}{c}\eta_{2}\\\eta_{4}\end{array}\right].
$$
Applying this operation recursively from the root to the leaves determines the incoming
trace on every leaf box.

Once $\phi_{\tau}$ is known on a leaf, its volume density is recovered independently by
\begin{equation}
\label{eq:recover_leaf_density}
q|_{\Omega_{\tau}}=Q^{\tau}\phi_{\tau}
\end{equation}
and the incoming field for box $\tau$ can be computed by
\begin{equation}
\label{eq:recover_leaf_incoming}
v|_{\Omega_{\tau}}=S_{\tau}\phi_{\tau}.
\end{equation}
The total field $u^{\rm tot} = u + u_{\rm in}$ inside the leaf box $\tau$ is recovered by superposing the
incoming field and the field induced by the volume density,
\begin{equation}
\label{eq:total_field_in_leaf}
u^{\rm tot}|_{\Omega_{\tau}}
=
\int_{\Omega_\tau} G_\kappa(\pxx,\pyy) q(\pyy) \, d\pyy  + v|_{\Omega_{\tau}}.
\end{equation}

Finally, the local scattered field is given by
\begin{equation}
\label{eq:scattered_field_in_leaf}
    u|_{\Omega_{\tau}} = u^{\rm tot}|_{\Omega_{\tau}} - u_{\rm in}|_{\Omega_{\tau}}.
\end{equation}

In the exterior of $\Omega$, the scattered field can be reconstructed
from the root outgoing trace
$$
\psi_{1}=F^{1}\phi_{1}
$$
by solving the radiating exterior Dirichlet problem.


The entire construction in Sections \ref{sec:local_scattering}--\ref{sec:thealgorithm}
is exact at the continuum level. Approximation enters only when the operators appearing
above are discretized, which is the subject of the next section.

\subsection{Storage efficiency}

There are several ways to reduce storage requirements. First, the scattering operators $F^\tau$ are needed only during the upward pass. Once a parent operator $F^\tau$ has been constructed, the scattering operators of its children may be discarded, since they are not used in the downward pass. Second, for recovery of the full interior field, the solution-transfer operators $T^\tau$ must be retained, while the leaf source operators $Q^\tau$ may either be stored or recomputed on demand during the downward pass. Recomputing the $Q^\tau$'s reduces memory usage at the expense of additional local solves. Finally, if only the root scattering operator $F^1$ is required, neither the $T^\tau$'s nor the $Q^\tau$'s need be retained; in that case all intermediate operators may be discarded as soon as they are no longer needed.

\section{Discretization and implementation details}
\label{sec:disc}

This section describes how we turn the continuum algorithm of
Section \ref{sec:thealgorithm} into an actual numerical method.

\subsection{Discretization}
\label{sec:tracediscretization}
All boundary traces such as $\{\phi_{\tau},\,\psi_{\tau},\eta,\gamma\}$
are represented via collocation. For a leaf box $\tau$, we place $p$ Gauss--Legendre
nodes on each side. Typically high $p$ is used, say $p=10$ or $p=20$.
For larger boxes, we use the composite Gauss--Legendre grids obtained by retaining
the boundary nodes of their children.

\begin{figure}
    \begin{minipage}{.5\linewidth}
    \centering
    \includegraphics[width=40mm]{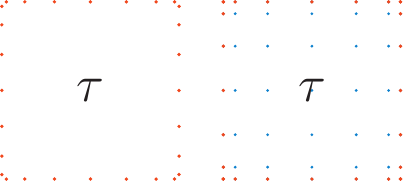}    
    \end{minipage}
    \hfill
    \begin{minipage}{80mm}
    \textit{The left picture shows the Gauss-Legendre nodes that we use to represent boundary traces, cf.~Section \ref{sec:tracediscretization}. The right shows the tensor product Chebyshev grid used for solving the local volume problem, cf.~Section \ref{sec:leaf_comp_continuum}.
    }
    \end{minipage}
    \caption{Local discretization on a leaf box.}
    \label{fig:placeholder}
\end{figure}

To facilitate the interior spectral solve, we also introduce a Chebyshev
discretization of $\Omega_{\tau}$. The tensor-product
Chebyshev grid provides the interior discretization on which the Lippman--Schwinger equation is solved. Local discretizations are shown in \cref{fig:placeholder}.

\subsection{Leaf computations}\label{subsec:leaf_comp}
\begin{figure}
    \centering
    \includegraphics[width=0.75\linewidth]{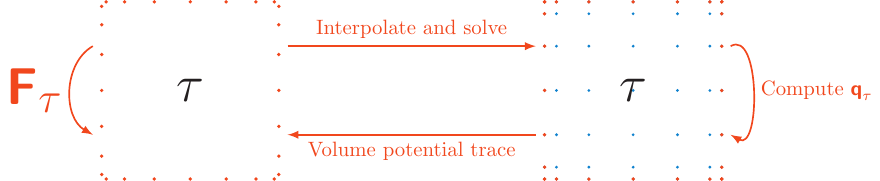}
    \caption{Decomposition of the discretized scattering operator $\mtx{F}_{\tau}$ into its constituent operators (see Section~\ref{subsec:leaf_comp}).}
    \label{fig:leaf_comp}
\end{figure}
As introduced in \cref{sec:leaf_comp_continuum}, the task is to find the scattering operator that maps the incoming field prescribed on the boundary of a leaf box $\tau$ to the corresponding outgoing scattered field on the same boundary.
The leaf computations are executed on tensor-product Chebyshev grids.

Recall that, in the continuous case, $F_{\tau}$ maps the incoming field $\phi_{\tau}$ to the outgoing field $\psi_{\tau}$, where $\psi_{\tau}$ is the boundary trace of $\int_{\Omega_{\tau}}G(\pxx,\pyy)q_{\tau}(\pyy)d\pyy$, and $q_{\tau}$ solves the local Lippmann--Schwinger equation with incoming field $\phi_{\tau}$. In discretized form, the computation of $\mtx{F}_{\tau}$ then involves the following steps:
\begin{itemize}
\item Interpolate the incoming trace from the Gauss--Legendre nodes to Chebyshev
boundary nodes, averaging the two side interpolants at each corner.
\item Use a local spectral collocation solve for the homogeneous Helmholtz equation
to reconstruct the incoming field at all tensor-product Chebyshev nodes.
\item Solve the local Lippmann--Schwinger equation using a high-order corrected
quadrature rule to obtain the (discretized) volume density $\vct{q}_{\tau}$.
\item Interpolate the trace of the volume potential generated by $\vct{q}_{\tau}$ to the Gauss--Legendre boundary nodes to obtain the outgoing trace.
\end{itemize}

Given an incoming field $\vct{\phi}_{\rm GL}$ on the Gauss--Legendre boundary
nodes, the operator $\mtx{S}$ reconstructs the corresponding field
$\vct{\phi}_{\rm vol}$ on the tensor-product Chebyshev grid. It consists of
interpolation from the Gauss--Legendre boundary nodes to the Chebyshev
boundary nodes, followed by a Chebyshev spectral collocation solve of the
homogeneous Helmholtz equation. Thus,
\begin{equation}
    \vct{\phi}_{\rm vol} = \mtx{S}\vct{\phi}_{\rm GL}.
\end{equation}

The induced volume density is then computed from the discretized
Lippmann--Schwinger equation,
\begin{equation}
    \mtx{L}\vct{q}_{\tau}
    =
    \mtx{B}\vct{\phi}_{\rm vol},
\end{equation}
where $\mtx{L} = \mtx{I} - \mtx{B} \mtx{G}$ is the discretized Lippmann--Schwinger operator, $\mtx{B}$ contains the scattering potential scaled by $-\kappa^2$  and
$\mtx{G}$ represents the volume Green's operator $q \mapsto \int_{\Omega_\tau} G_\kappa(\cdot, \pyy) q(\pyy) \, d\pyy$ discretized using
a high-order corrected quadrature rule. Hence,
\begin{equation}
    \vct{q}_{\tau}
    =
    \mtx{L}^{-1}
    \mtx{B}\mtx{S}\vct{\phi}_{\rm GL}.
\end{equation}

Finally, the scattered field generated by $q$ is evaluated on the
Chebyshev boundary nodes and interpolated to the Gauss--Legendre boundary
nodes. Denoting the resulting operator by $\mtx{M}$, we have
\begin{equation}
    \vct{\psi}_{\rm GL} = \mtx{M}\vct{q}_{\tau}.
\end{equation}
Combining the three stages gives
\begin{equation}
    \vct{\psi}_{\rm GL}
    =
    \mtx{M}
    \mtx{L}^{-1}
    \mtx{B}\mtx{S}\vct{\phi}_{\rm GL}.
\end{equation}
Therefore, the scattering matrix for leaf $\tau$ is
\begin{equation}
    \boxed{
    \mtx{F}^{\tau}
    =
    \mtx{M}
    \mtx{L}^{-1}
    \mtx{B}\mtx{S}
    }.
\end{equation}

The leaf continuation operators $\mtx{S}$, scattering trace operator $\mtx{M}$ and the high-order quadrature matrices $\mtx{G}$ are
\textit{\textbf{universal}} when all leaf boxes have the same size. They therefore need to be constructed only once. 

\subsection{Discretization of the merge stage}\label{sec:disc_merge}
As already noted in Remark~\ref{rem:universal_ops}, the operators $E_{\tau},U_{\tau}$ and $R_{\tau}$ from~\eqref{eq:generic_merge} are universal at each level \textit{where they exist}. At level $\ell = L$, there is no continuation operator (leaf nodes have no children), so we can index all discretized continuation operators by $\mtx{U}^{(\ell)}, \ell=0,\ldots,L-1$. Similarly (the root node has no siblings) we can index all sibling-evaluation operators by $\mtx{E}^{(\ell)},\mtx{R}^{(\ell)},\ell=1,\ldots,L$. For notational convenience we will drop the index $\ell$ for the blocks and operators that constitute $\mtx{E}^{(\ell)},\mtx{R}^{(\ell)}$ and $\mtx{U}^{(\ell)}$.

\begin{figure}
\begin{minipage}{.3\linewidth}
    \includegraphics[width=\linewidth]{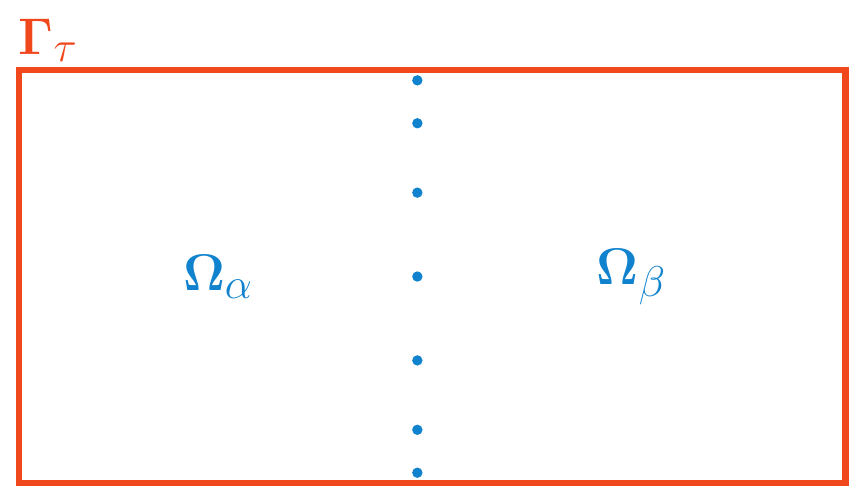}\\
    \includegraphics[width=\linewidth]{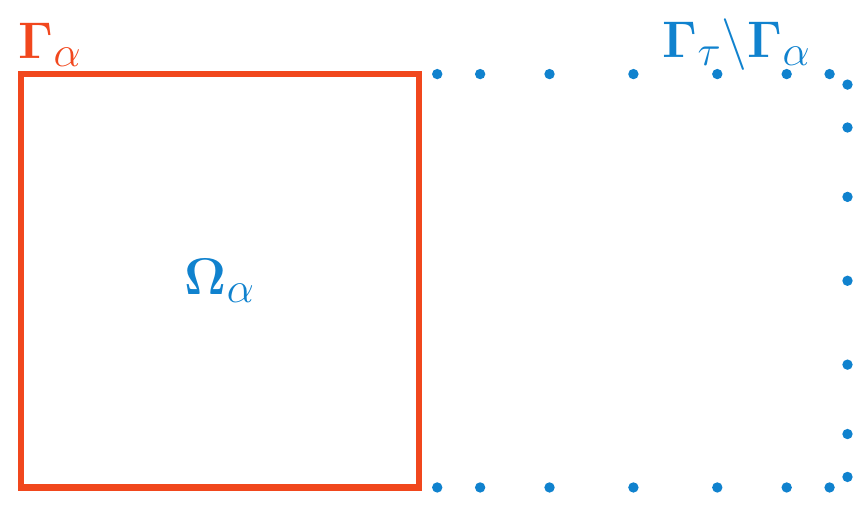}
\end{minipage}
\hfill
\begin{minipage}{.6\linewidth}
    \textit{Continuation and sibling-evaluation for nodes $\alpha,\beta$ at level $\ell$. Displayed are the set-up for $\mtx{U}^{(\ell-1)}$ (top) and the blocks $\mtx{E}_{21},\mtx{E}_{23}$ in $\mtx{E}^{(\ell)}$ (bottom). In both cases a constant-coefficient Helmholtz BVP is discretized and solved using a BIE method on the solid red boundary. For $\mtx{U}^{(\ell-1)}$ the BVP is an interior Helmholtz equation with Dirichlet boundary conditions. For $\mtx{E}_{21},\mtx{E}_{23}$ the BVP is an exterior Helmholtz equation with Dirichlet conditions and Sommerfeld radiation condition. The evaluation points are indicated in blue. For spectral discretization, the Gauss-Legendre points approach the boundary linearly as the polynomial order $p$ increases forcing the usage of BIE solvers with close-evaluation capabilities.}
\end{minipage}
    \caption{Set-up for the discretization of $\mtx{U}^{(\ell-1)}$ and $\mtx{E}^{(\ell)}$.}
    \label{fig:discEandU}
\end{figure}

As can be seen from~\eqref{eq:E_adjacent} and~\eqref{eq:R_structure}, $\mtx{E}^{(\ell)}$ and $\mtx{R}^{(\ell)}$ can be written as
\begin{equation}
\label{eq:E_ell}
\mtx{E}^{(\ell)}=
\left[\begin{array}{cccc}
\mtx{0}&\mtx{E}_{12}&\mtx{0}&\mtx{E}_{14}\\
\mtx{E}_{21}&\mtx{0}&\mtx{E}_{23}&\mtx{0}\\
\mtx{0}&\mtx{0}&\mtx{0}&\mtx{I}\\
\mtx{0}&\mtx{0}&\mtx{I}&\mtx{0}
\end{array}\right],\quad
\mtx{R}^{(\ell)}=
\left[\begin{array}{cccc}
\mtx{I}&\mtx{E}_{12}&\mtx{0}&\mtx{E}_{14}\\
\mtx{E}_{21}&\mtx{I}&\mtx{E}_{23}&\mtx{0}
\end{array}\right].
\end{equation}
Figure~\ref{fig:discEandU} shows the set-up for $\mtx{U}^{(\ell-1)}$ and $\mtx{E}_{21},\mtx{E}_{23}$. The construction for $\mtx{E}_{12},\mtx{E}_{14}$ is similar.
Computationally, the structure of $[\mtx{E}_{21},\mtx{E}_{23}]$ is
$$[\mtx{E}_{21},\mtx{E}_{23}] = \mtx{BIE}^{(\rm ext)}_{\rm eval}\circ\mtx{BIE}^{(\rm ext)}_{\rm solve}\circ\mtx{Interp}.$$
Here $\mtx{Interp}$ maps the boundary data, collocated at piecewise Gauss-Legendre nodes, to interpolated boundary data at a point distribution on the boundary that is suited for a BIE solver. The exterior BIE solver $\mtx{BIE}^{(\rm ext)}_{\rm solve}$ solves for a (discrete) surface density, whose precise form depends on the chosen formulation. Finally, $\mtx{BIE}^{(\rm ext)}_{\rm eval}$ maps the resulting surface density to the values of the Helmholtz solution at the specified exterior target DOFs. The matrix $\mtx{U}^{(\ell-1)}$ has the same structure, but interior:
$$\mtx{U}^{(\ell-1)} = \mtx{BIE}^{(\rm int)}_{\rm eval}\circ\mtx{BIE}^{(\rm int)}_{\rm solve}\circ\mtx{Interp}.$$
In our implementation we use Brakhage–Werner-type indirect BIE formulations for both the exterior and interior problem. Concretely:
\begin{equation*}
\mtx{BIE}^{(\rm int)}_{\rm solve}=\left(\mtx{DL}+i\xi\mtx{SL}+\frac{1}{2}\mtx{I}\right)^{-1},\quad
    \mtx{BIE}^{(\rm ext)}_{\rm solve}=\left(\mtx{DL}+i\xi\mtx{SL}-\frac{1}{2}\mtx{I}\right)^{-1}    
\end{equation*} 
with $\mtx{DL}$ and $\mtx{SL}$ discretizations of the boundary double-layer and boundary single-layer operators (here $\mtx{SL}$ is not to be confused with the product of $\mtx{S}$ and $\mtx{L}$ from Section~\ref{subsec:leaf_comp}) and $\xi\in\mathbb{R}$ is a scalar usually chosen to be $\pm \kappa$. We construct and solve both using the chunkIE package (see~\cite{chunkie}).

For a fixed level $\ell$, $\mtx{U}^{(\ell)}$ and $\mtx{E}^{(\ell)}$ are defined on the same boundary. As such, the matrices $\mtx{DL}$ and $\mtx{SL}$ can be re-used for the construction of $\mtx{E}^{(\ell)}$ and $\mtx{U}^{(\ell)}$.
\subsection{Complexity}\label{subsec:complexity}
The complexity of our proposed solver is dominated by the following three terms:
\begin{enumerate}
    \item $c_{\rm leaf}$: the construction, at the leaf level, of $\mtx{Q}_{\tau}$ and $\mtx{F}_{\tau}$. 
    \item $c_{\rm pre}$: the off-line computation, at each level, of $\mtx{E}^{(\ell)},\mtx{U}^{(\ell)}$ and $\mtx{R}^{(\ell)}$.
    \item $c_{\rm inv}$: the inversion, at each level, of $(\mtx{I}-\mtx{E}^{(\ell)}\mtx{F})$.
\end{enumerate}
We set $N=p^2\cdot2^{L}$ to be the total number of \textit{unreduced} degrees of freedom. For simplicity we assume that $L$ is even. 

The leaf scattering computation $c_{\rm leaf}$ is dominated by the inversion of $\mtx{L}_{\tau}$, which is of order $\mathcal{O}(p^6)$. The total cost $c_{\rm leaf}$ is therefore 
$$c_{\rm leaf} = \mathcal{O}(p^62^{L}) = \mathcal{O}(p^4N).$$
Usually, the leaf level matrices are quite small. Even at $p=30$, they are only of size $900 \times 900$, making them perfectly suited for highly efficient BLAS level $3$ operations.

The cost $c_{\rm inv}$, by a geometric series argument, is bounded by a constant times the cost of the inversion of $(\mtx{I}-\mtx{E}^{(0)}\mtx{F})$ at the root, where the size of the system is precisely $2^{L/2}p \times 2^{L/2} p = \sqrt{N}\times\sqrt{N}$. Then, trivially,
$$c_{\rm inv} = \mathcal{O}(N^{3/2}).$$

The cost $c_{\rm pre}$ depends heavily on the chosen BIE solvers. We present the argument for the continuation matrices $\mtx{U}^{(\ell)}, \ell=1,\ldots,L$ here, but the same argument applies to the sibling evaluation operators.

Assume that the system matrices $\mtx{BIE}_{\rm solve}^{\rm (int)}$ at level $\ell$ are square and of size $n_{\rm BIE}^{(\ell)}\times n_{\rm BIE}^{(\ell)}$. With $c_{\ell}$ the number of boundary segments of a box at level $\ell$, we ensure $n_{\rm BIE}^{(\ell)} = \mathcal{O}(c_{\ell}p) = \mathcal{O}(2^{(L-\ell)/2}p)$. Naive construction and inversion of the BIE system results in a cost of $\mathcal{O}((n_{\rm BIE}^{(\ell)})^2+(n_{\rm BIE}^{(\ell)})^3)$. The cost of applying interpolation and evaluation by $\mtx{Interp},\mtx{BIE}^{\rm int}_{\rm eval}$ scales as follows: at level $\ell$ these operators are block diagonal matrices, consisting of $\mathcal{O}(2^{L-\ell})$ blocks. Each of the blocks in the interpolation (resp. evaluation) can be of variable row size (resp. column size), but always has column size (resp. row size) $p$. This means the cost multiplying $\mtx{BIE}_{\rm solve}^{\rm(int)}$ left and right by these block diagonal matrices is bounded by $\mathcal{O}(p (n_{\rm BIE}^{(\ell)})^2)$.

Summing over levels, again using a geometric series argument, finally gives us
$$c_{\rm pre} = \mathcal{O}(\sum_{\ell}(p(n_{\rm BIE}^{(\ell)})^2+(n_{\rm BIE}^{(\ell)})^3)) = \mathcal{O}(N^{3/2}).$$
The final estimate on the right-hand side obfuscates the dependence of the off-line phase on the chosen BIE solver, and for this reason we present both complexity estimates. While we have presented the asymptotic cost for naive BIE construction and inversion, at large scale both can be accelerated using fast algorithms like the Fast Multipole Method (FMM, see, e.g.,~\cite{rokhlin1987,rokhlin_1990,ying2004kernel}), hierarchical matrix techniques (see~\cite{2008_bebendorf_book,2003_hackbusch,2005_martinsson_fastdirect}) and recursive compression (see~\cite{RCIP}).

\section{Numerical examples}
\label{sec:num}

We present six numerical experiments with various wavenumbers $\kappa$ and scattering potentials $b(\pxx)$. The computational domain is the unit square $[0,1]^2$. The numerical computations in this section were carried out using MATLAB. All timings are reported on a workstation with an Intel(R) Xeon(R) Gold 6248R processor (48 cores with 2 threads per core at 3.00\,GHz base frequency) and 1.5\,TB of RAM. Results are presented in \cref{tab:result_time,tab:result_error}.

We note the following quantities:
    \begin{center}
    \begin{tabular}{@{} l  l @{}}
        $n_{\rm level}$     & The number of levels of the hierarchical tree. \\
        $n_\text{leaves}$   & The number of leaves on the finest level. \\
        $N_\text{tot}$      & The total number of degrees of freedom. \\
        $t_\text{precomp}$  & Time in seconds to precompute the $\mtx{S}, \mtx{M},\mtx{G},\mtx{E}, \mtx{U},\mtx{R}$ matrices. \\
        $t_\text{build}$    & Time in seconds to build scattering matrices $\mtx{F}$. \\
        $t_\text{eval}$    & Time in seconds to evaluate exterior field. \\
        $t_\text{solve}$    & Time in seconds to solve for interior field. \\
        $E_\text{ext}$      & The relative error in the exterior field (total field) evaluated at selected points.\\
        $E_\text{int}$      & The relative error in the interior field (total field) evaluated at selected off-grid points.\\
        $E_\text{grid}$     & The maximal relative error in the interior field (total field) among all tensor Chebyshev grid points.\\
    \end{tabular}
    \end{center}
    All errors are in the $L^{\infty}$-norm and obtained by comparing to a reference solution at $n_{\rm level} = 17$ and $p=20$. The total field refers to the sum of the scattered and incoming fields, i.e.,  $u^{\rm tot} = u + u_{\rm in}$, as in section~\ref{sec:downward-pass}.

    We consider six test problems, where the incident field is given by a plane wave $\puu_\text{in}(\pxx) = \exp(\text{i}\kappa (x_1-0.5))$:
    \begin{center}
    \begin{tabular}{@{}p{0.20\textwidth}p{0.75\textwidth}@{}}

\textit{Gaussian bump:}
&
A Gaussian bump given by
$b(\pxx) = 1.5\exp(-160 \lVert \pxx - \pxx_{c} \rVert^2)$
with $\kappa=25$, where $\pxx_c$ is the center of the computational domain (\cref{fig:bfun}(a)). The computational domain is about $4$ wavelengths across.
\\[1em]

\textit{Smooth mollifier:}
&
A smooth mollifier given by
$b(\pxx) = \exp(-\frac{\epsilon^2}{\epsilon^2-\lVert \pxx - \pxx_{c} \rVert^2})$ for $\lVert \pxx - \pxx_{c} \rVert < \epsilon$ with $ \epsilon = 0.4$ (\cref{fig:bfun}(b)). The wavenumber is set to $\kappa = 25$, so that the computational domain is about $4$ wavelengths across.
\\[1em]

\textit{Cavity:}
&
A cavity potential is given in polar coordinates by $b(r,\theta) = (1-\sin(0.5\theta)^{500})\exp(-2000(0.1-r^2)^2)$ (\cref{fig:bfun}(c)). The wavenumber is set to $\kappa = 50.27$, so that the computational domain is about $8$ wavelengths across.
\\[1em]

\textit{Lens:}
&
A vertically-graded lens given by
$b(\pxx) = 4(x_2-0.7)\left[1-\operatorname{erf}\left(25(\|\pxx - \pxx_c\|-0.3)\right)\right]$
with $\kappa=200$ (\cref{fig:bfun}(d)). The computational domain is about $31.8$ wavelengths across.
\\[1em]

\textit{Photonic crystal:}
&
A 20 by 20 square grid of Gaussian bumps with a channel removed
with $\kappa=77.1$ (the first complete band gap) (\cref{fig:bfun}(e)). The computational domain is about $12.2$ wavelengths across.
\\[1em]

\textit{Pyramid:}
&
A pyramid potential given by $b(\pxx) = \max\bigl(1-\frac{|x_1-0.5|+|x_2-0.5|}{0.4},0\bigr)$ with $\kappa=25$ (\cref{fig:bfun}(f)). The computational domain is about $4$ wavelengths across.
\\

\end{tabular}
\end{center}

\cref{tab:result_time} reports the computational cost of the algorithm for increasing numbers of levels. As the hierarchy is refined from 3 to 15 levels, the precomputation time increases from 13.5 s to 49.3 s, reflecting the additional cost associated with constructing the operators required at each level of the hierarchy. In contrast, the evaluation time remains nearly constant despite the substantial increase in problem size. The build and solve times exhibit more pronounced growth, roughly linear in terms of $N_{\rm tot}$. These results demonstrate that the computational cost of the hierarchical merge is dominated increasingly by the construction and solution stages as the number of levels increases, whereas the evaluation stage remains relatively insensitive to the overall problem size.

\cref{tab:result_error} summarizes the convergence behavior of the hierarchical merge for six scattering potentials. For the Gaussian bump, smooth mollifier, and cavity test cases, the exterior, interior on-grid and off-grid errors decrease rapidly with refinement, reaching errors between approximately $10^{-10}$ and $10^{-8}$ at the finest resolution. The lens and photonic-crystal examples are more demanding on coarse grids and require substantially finer discretizations before achieving high accuracy. In particular, the lens errors decrease from $O(1)$ at three levels to below $10^{-7}$ at fifteen levels, while the photonic-crystal errors exhibit a substantial reduction at eleven levels, and decrease gradually until the finest resolution. Additionally the BIE operators from Section~\ref{sec:disc_merge} are refined in this case (i.e., $n_{\rm BIE}$ is increased, see Section~\ref{subsec:complexity}). The pyramid case, since its $b$-function is only $C^0$, converges more gradually. Overall, these results demonstrate that the hierarchical merge produces convergent solutions across a range of smooth and highly structured scattering potentials, while the resolution required to attain a given accuracy depends on the complexity of the scattering potential.
For the smooth examples, the errors eventually plateau under further refinement. We suspect that this behavior is primarily due to limitations in the accuracy of the boundary integral evaluations used to construct the universal operators, particularly for targets very close to the boundary.

\begin{table}[htp]
  \centering\small
  \caption{Computational cost for different numbers of levels, at $\kappa=77.1$ and $p=20$.}
  \label{tab:result_time}
  \begin{minipage}{.9\linewidth}
  \renewcommand{\footnoterule}{}
  \centering
  \sisetup{
    table-format=1.2e-1,
    table-auto-round=true,
    tight-spacing=true
  }

  \begin{tabular}{ccc|S[retain-zero-exponent=true]S[retain-zero-exponent=true]S[retain-zero-exponent=true]S[retain-zero-exponent=true]}
    \toprule
    \multicolumn{3}{c|}{\bf Number of} &
    \multicolumn{4}{c}{\bf Time [s]} \\
    levels & leaves & points &
    {precomp.} & {build} & {eval.} & {solve} \\
    $n_{\rm level}$ & $n_\text{leaves}$ & $N_\text{tot}$ &
    {$t_\text{precomp}$} & {$t_\text{build}$} &
    {$t_\text{eval}$} & {$t_\text{solve}$} \\
    \midrule

    \num{3} & \num{4} & \num{1600}
    & \num{1.35e+01} & \num{5.80e-02} & \num{4.16e-2} & \num{5.11e-03} \\
    \num{5} & \num{16} & \num{6400}
    & \num{1.75e+01} & \num{2.49e-01} & \num{4.36e-2} & \num{1.09e-02} \\
    \num{7} & \num{64} & \num{25600}
    & \num{2.14e+01} & \num{9.18e-01} & \num{4.58e-2} & \num{4.26e-02} \\
    \num{9} & \num{256} & \num{102400}
    & \num{2.59e+01} & \num{4.70e+00} & \num{3.88e-2} & \num{1.82e-01} \\
    \num{11} & \num{1024} & \num{409600}
    & \num{3.04e+01} & \num{1.96e+01} & \num{4.14e-2} & \num{6.80e-01} \\
    \num{13} & \num{4096} & \num{1638400}
    & \num{3.78e+01} & \num{8.36e+01} & \num{5.10e-2} & \num{2.48e+00} \\
    \num{15} & \num{16384} & \num{6553600}
    & \num{4.93e+01} & \num{3.74e+02} & \num{7.05e-2} & {\num{2.18e+01}\footnote{At this scale, $\sim 45$GB is needed to store the full system. The jump in solve cost is due to system-level memory effects.}} \\
    \bottomrule
  \end{tabular}
\end{minipage}
\end{table}
\clearpage
\begin{figure}
    \centering
    \subfigure[Gaussian bump]{
    \includegraphics[width=.45\linewidth,trim=1cm 6cm 2cm 7cm,clip]{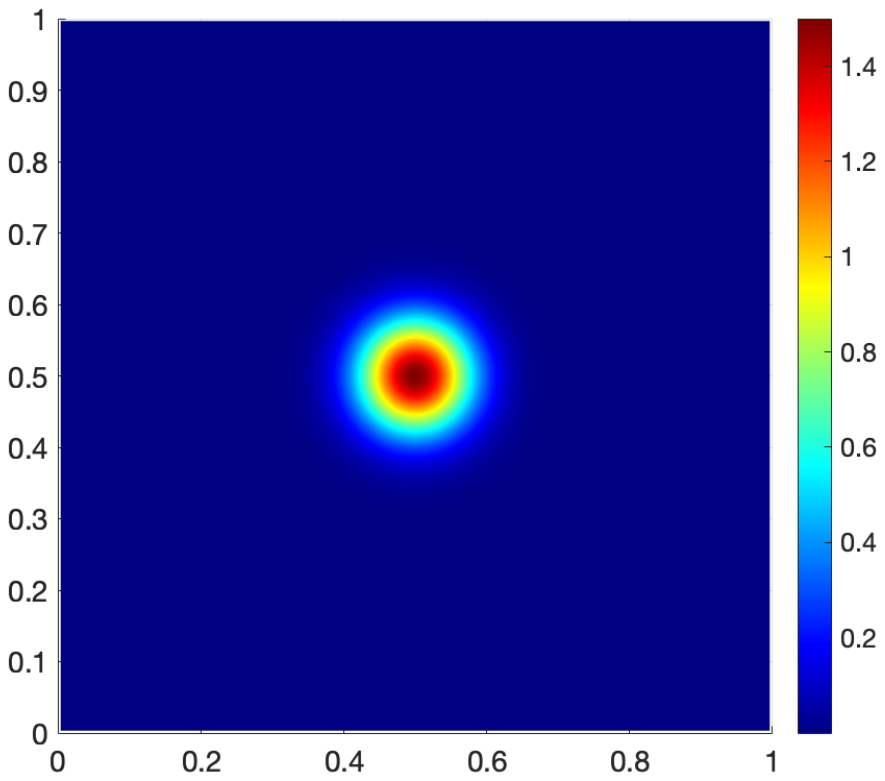}
    \label{fig:bfun_gaussianbump}
    } %
    \subfigure[Smooth mollifier]{
    \includegraphics[width=.45\linewidth,trim=1cm 6cm 2cm 7cm,clip]{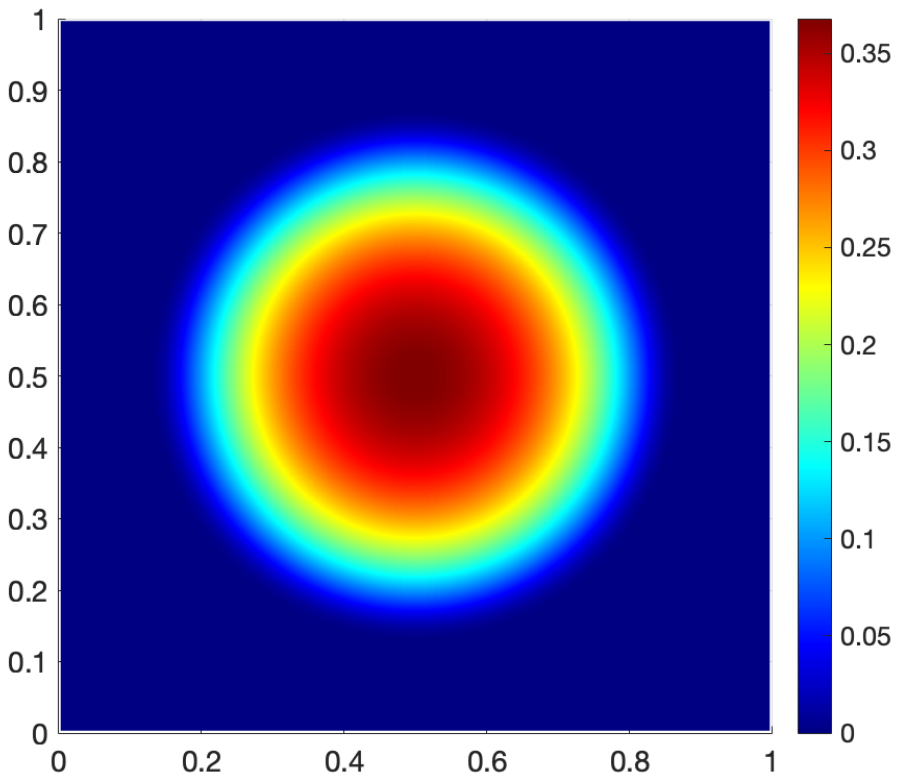}
    \label{fig:bfun_smoothmollifier}
    }
    \bigskip
    \subfigure[Cavity]{
    \includegraphics[width=.45\linewidth,trim=1cm 6cm 2cm 7cm,clip]{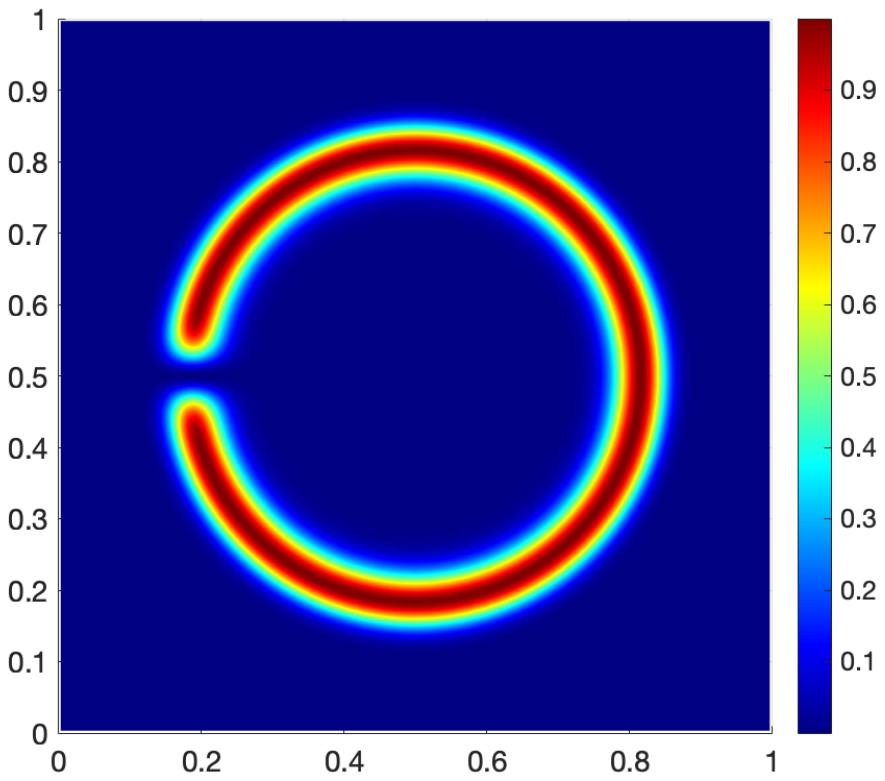}
    \label{fig:bfun_cavity}
    } %
    \subfigure[Lens]{
    \includegraphics[width=.45\linewidth,trim=1cm 6cm 2cm 7cm,clip]{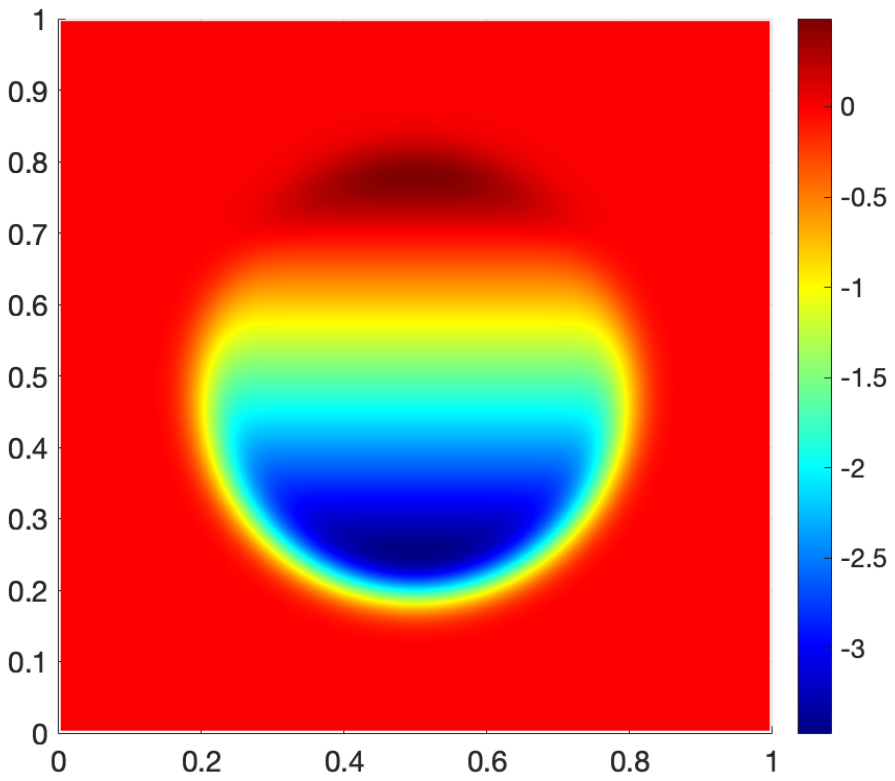}
    \label{fig:bfun_lens}
    }
    \bigskip
    \subfigure[Photonic crystal]{
    \includegraphics[width=.45\linewidth,trim=1cm 6cm 2cm 7cm,clip]{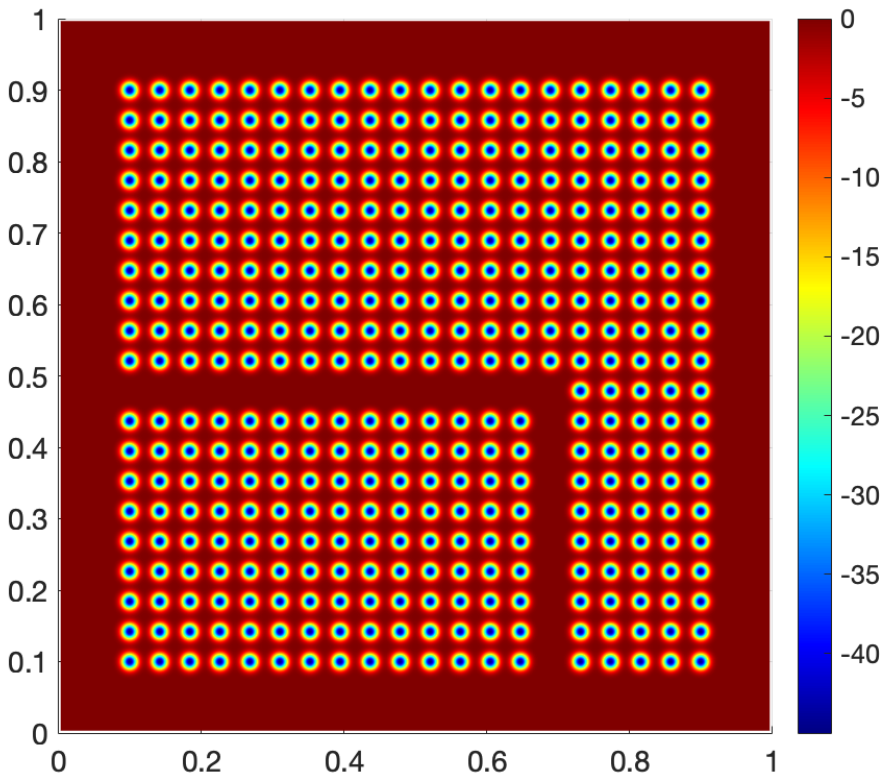}
    \label{fig:bfun_photonic_crystal}
    } %
    \subfigure[Pyramid]{
    \includegraphics[width=.45\linewidth,trim=1cm 6cm 2cm 7cm,clip]{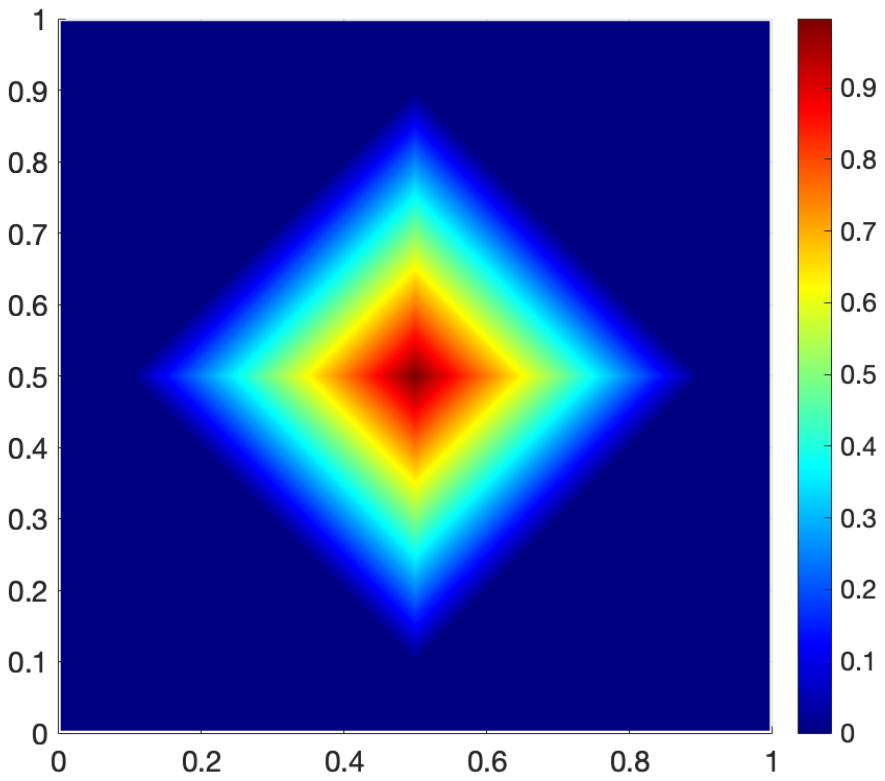}
    \label{fig:bfun_pyramid}
    }
    \bigskip
    \caption{Scattering potentials $b(\pxx)$.}
    \label{fig:bfun}
\end{figure}

\begin{figure}
    \centering
    \subfigure[Gaussian bump]{
    \includegraphics[width=.45\linewidth,trim=1cm 6cm 2cm 7cm,clip]{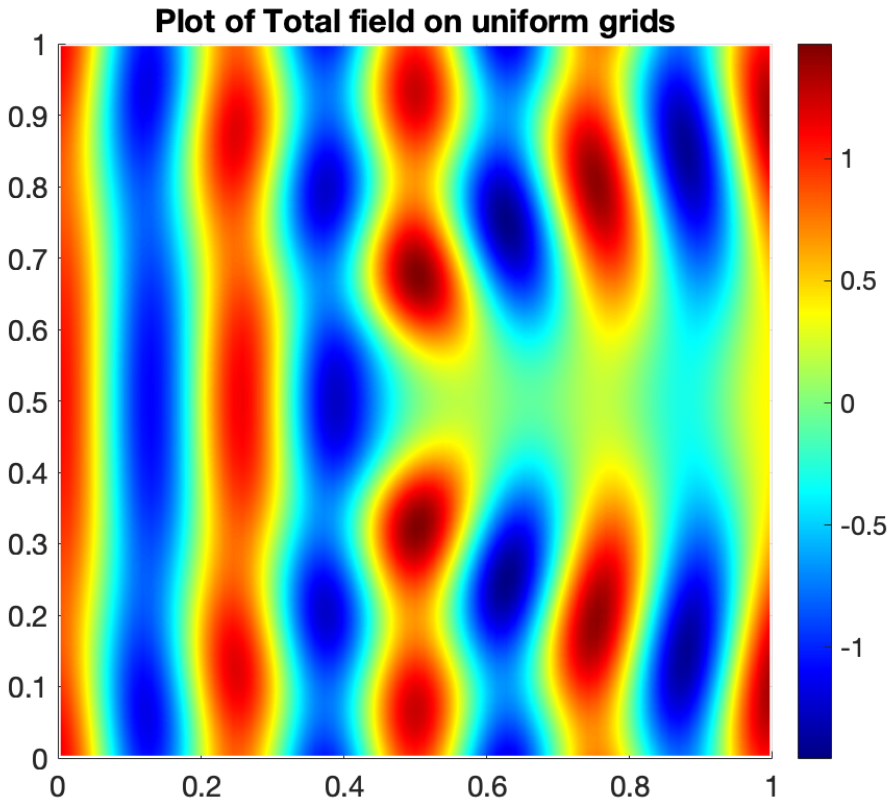}
    } %
    \subfigure[Smooth mollifier]{
    \includegraphics[width=.45\linewidth,trim=1cm 6cm 2cm 7cm,clip]{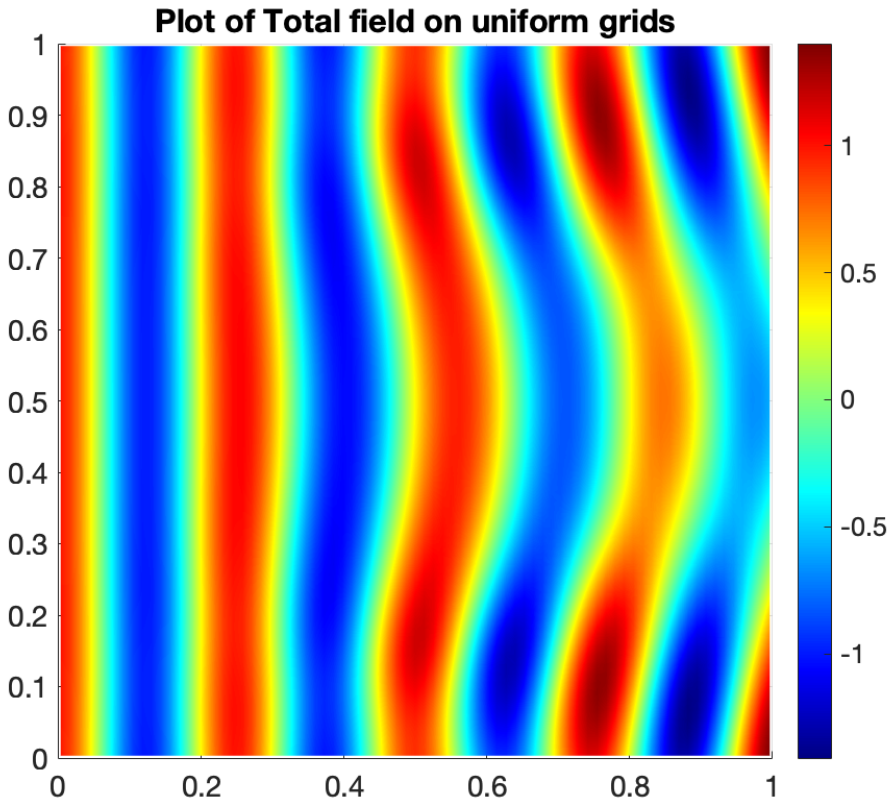}
    }
    \bigskip
    \subfigure[Cavity]{
    \includegraphics[width=.45\linewidth,trim=1cm 6cm 2cm 7cm,clip]{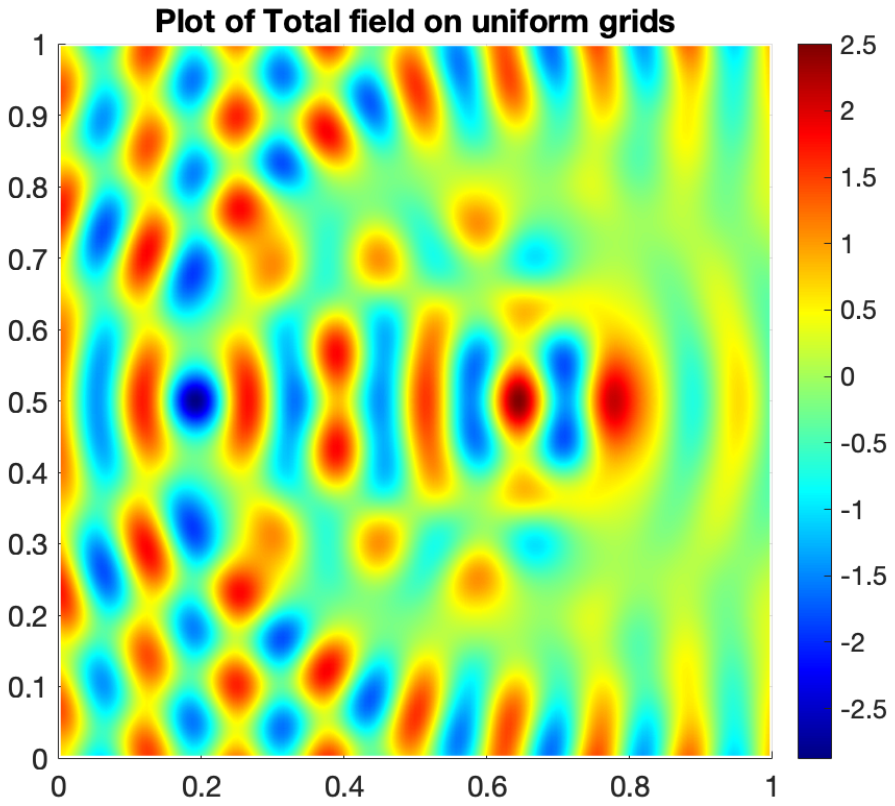}
    } %
    \subfigure[Lens]{
    \includegraphics[width=.45\linewidth,trim=1cm 6cm 2cm 7cm,clip]{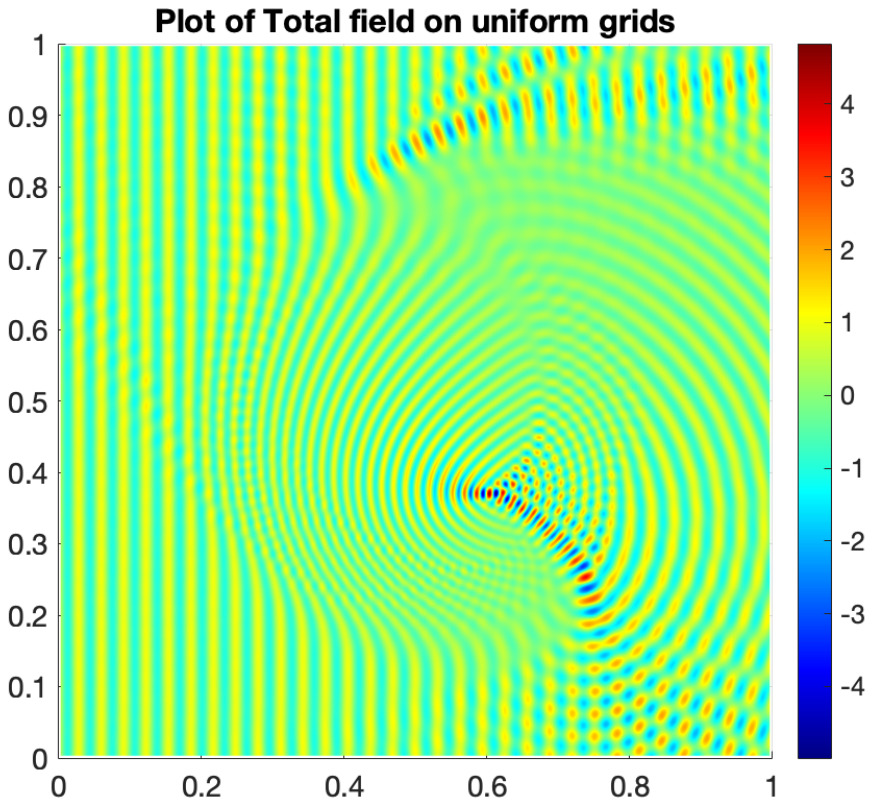}
    }
    \bigskip
    \subfigure[Photonic crystal]{
    \includegraphics[width=.45\linewidth,trim=1cm 6cm 2cm 7cm,clip]{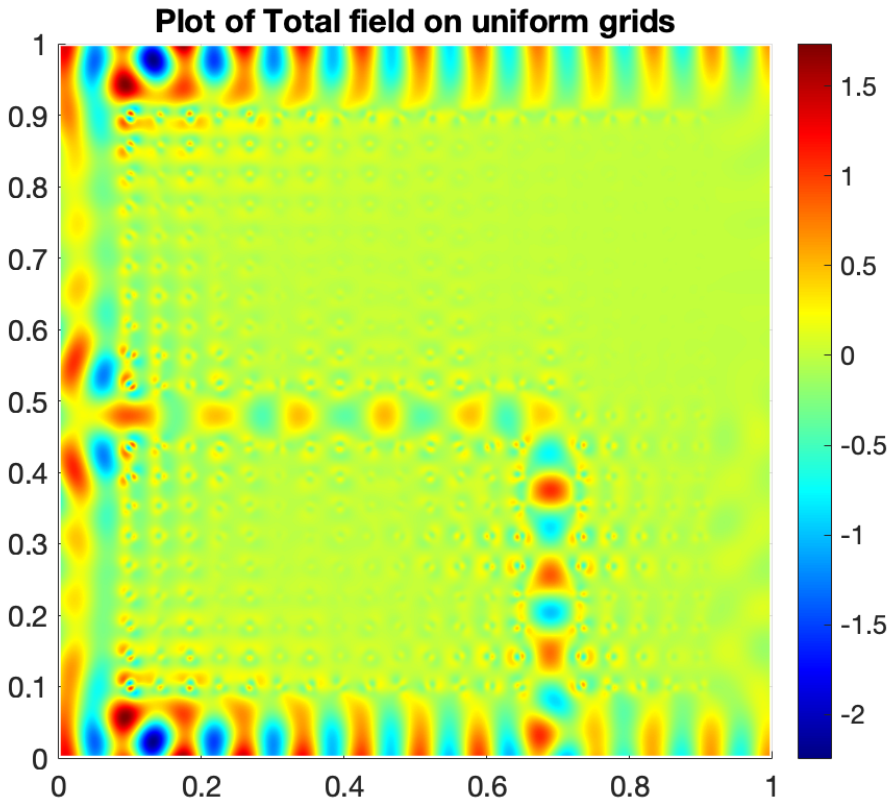}
    } %
    \subfigure[Pyramid]{
    \includegraphics[width=.45\linewidth,trim=1cm 6cm 2cm 7cm,clip]{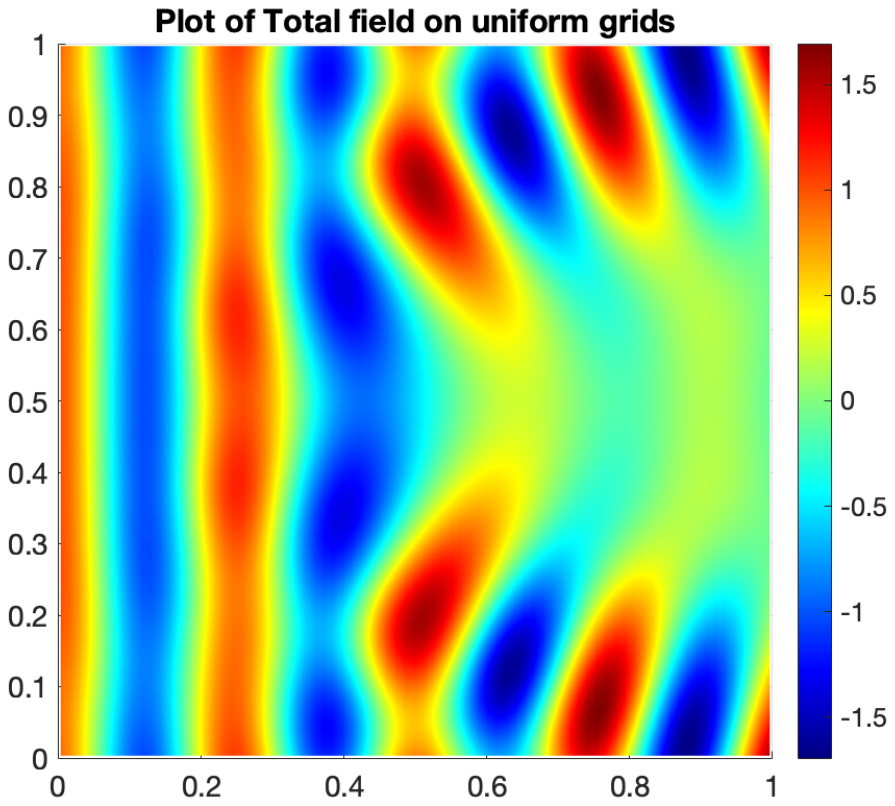}
    }
    \bigskip
    \caption{Total field for different scattering potentials $b(\pxx)$, given the incoming plane wave $\puu_\text{in}(\pxx) = \exp(\text{i}\kappa (x_1-0.5))$.}
    \label{fig:utot}
\end{figure}
\clearpage

\begin{table}[htp]
  \centering\small
  \caption{Results of the proposed method for six scattering potentials.}
  \label{tab:result_error}
  \vspace*{5pt}

  \sisetup{
    table-format=1.2e-1,
    table-auto-round=true,
    tight-spacing=true
  }
\begin{minipage}{.8\linewidth}
\renewcommand{\footnoterule}{}
\centering
  \begin{tabular}{l|ccc|S[retain-zero-exponent=true]S[retain-zero-exponent=true]S[retain-zero-exponent=true]}
    \toprule
    \multicolumn{1}{c|}{\bf Scattering potential} &
    \multicolumn{3}{c|}{\bf Number of} &
    \multicolumn{3}{c}{\bf Error} \\
    & levels & leaves & points &
    {ext.} & {int.} & {grid} \\
    & $n_{\rm level}$ & $n_{\text{leaves}}$ & $N_{\text{tot}}$ &
    {$E_{\text{ext}}$} & {$E_{\text{int}}$} & {$E_{\text{grid}}$} \\
    \midrule

    \multirow{5}{*}{Gaussian bump}
      & \num{3}  & \num{4}    & \num{1600}   & \num{3.11e-07} & \num{3.54e-05} & \num{1.89e-05} \\
      & \num{5}  & \num{16}   & \num{6400}   & \num{1.22e-09} & \num{2.47e-07} & \num{5.55e-07} \\
      & \num{7}  & \num{64}   & \num{25600}  & \num{2.00e-10} & \num{2.33e-08} & \num{1.58e-07} \\
      & \num{9}  & \num{256}  & \num{102400} & \num{1.73e-10} & \num{3.86e-08} & \num{8.22e-08} \\
      & \num{11} & \num{1024} & \num{409600} & \num{1.88e-10} & \num{1.65e-08} & \num{3.53e-08} \\

    \midrule

    \multirow{5}{*}{Smooth mollifier}
      & \num{3}  & \num{4}    & \num{1600}   & \num{5.69e-05} & \num{2.90e-04} & \num{2.83e-04} \\
      & \num{5}  & \num{16}   & \num{6400}   & \num{1.41e-06} & \num{8.88e-06} & \num{9.69e-06} \\
      & \num{7}  & \num{64}   & \num{25600}  & \num{3.60e-10} & \num{8.03e-08} & \num{1.06e-07} \\
      & \num{9}  & \num{256}  & \num{102400} & \num{3.94e-10} & \num{2.43e-08} & \num{3.36e-08} \\
      & \num{11} & \num{1024} & \num{409600} & \num{4.26e-10} & \num{7.23e-09} & \num{3.47e-08} \\

    \midrule

    \multirow{5}{*}{Cavity}
      & \num{3}  & \num{4}    & \num{1600}   & \num{2.50e-02} & \num{6.76e-02} & \num{6.24e-02} \\
      & \num{5}  & \num{16}   & \num{6400}   & \num{3.55e-05} & \num{2.90e-04} & \num{1.33e-04} \\
      & \num{7}  & \num{64}   & \num{25600}  & \num{1.04e-09} & \num{9.91e-08} & \num{3.94e-07} \\
      & \num{9}  & \num{256}  & \num{102400} & \num{8.23e-10} & \num{3.08e-08} & \num{2.90e-07} \\
      & \num{11} & \num{1024} & \num{409600} & \num{7.49e-10} & \num{4.14e-09} & \num{1.61e-07} \\

    \midrule

    \multirow{7}{*}{Lens\footnote{At $\kappa = 200, p=20$, we need $16\times 16$ leaf boxes to ensure $10$ DOFs per wave length.}}
    & \num{3} & \num{4} & \num{1600}
    & \num{1.00e+01} & \num{1.75e+01} & \num{3.35e+01} \\
    & \num{5} & \num{16} & \num{6400}
    & \num{1.01e+01} & \num{8.59e+02} & \num{2.54e+03} \\
    & \num{7} & \num{64} & \num{25600}
    & \num{5.19e-01} & \num{9.52e-01} & \num{3.67e+00} \\
    & \num{9} & \num{256} & \num{102400}
    & \num{2.37e-05} & \num{4.62e-04} & \num{3.10e-04} \\
    & \num{11} & \num{1024} & \num{409600}
    & \num{1.36e-07} & \num{8.45e-07} & \num{3.87e-06} \\
    & \num{13} & \num{4096} & \num{1638400}
    & \num{6.16e-08} & \num{1.76e-07} & \num{1.34e-06} \\
    & \num{15} & \num{16384} & \num{6553600}
    & \num{4.52e-08} & \num{1.12e-07} & \num{9.03e-07} \\

    \midrule

    \multirow{7}{*}{Photonic crystal\footnote{Due to the strong localization of the Gaussian bumps, at $p=20$, $32\times 32$ leaf boxes are needed to resolve $b(\pxx)$.}}
    & \num{3} & \num{4} & \num{1600}
    & \num{1.09e+00} & \num{2.54e+00} & \num{5.98e+00} \\
    & \num{5} & \num{16} & \num{6400}
    & \num{4.09e-01} & \num{8.49e-01} & \num{1.43e+00} \\
    & \num{7} & \num{64} & \num{25600}
    & \num{3.45e-01} & \num{1.27e+00} & \num{2.68e+01} \\
    & \num{9} & \num{256} & \num{102400}
    & \num{1.26e-01} & \num{3.09e-01} & \num{6.25e+00} \\
    & \num{11} & \num{1024} & \num{409600}
    & \num{1.17e-06} & \num{4.33e-05} & \num{3.08e-05} \\
    & \num{13} & \num{4096} & \num{1638400}
    & \num{3.48e-08} & \num{2.25e-07} & \num{2.26e-06} \\
    & \num{15} & \num{16384} & \num{6553600}
    & \num{1.05e-08} & \num{4.68e-08} & \num{4.57e-07} \\

    \midrule

    \multirow{7}{*}{Pyramid\footnote{The pyramid case converges more gradually, as it has limited regularity, which is only $C^0$.}}
    & \num{3} & \num{4} & \num{1600}
    & \num{1.08e-03} & \num{2.35e-03} & \num{2.26e-03} \\
    & \num{5} & \num{16} & \num{6400}
    & \num{2.38e-04} & \num{6.89e-04} & \num{1.00e-03} \\
    & \num{7} & \num{64} & \num{25600}
    & \num{4.84e-05} & \num{1.23e-04} & \num{4.17e-04} \\
    & \num{9} & \num{256} & \num{102400}
    & \num{1.23e-05} & \num{2.32e-05} & \num{1.19e-04} \\
    & \num{11} & \num{1024} & \num{409600}
    & \num{2.57e-06} & \num{5.07e-06} & \num{2.76e-05} \\
    & \num{13} & \num{4096} & \num{1638400}
    & \num{7.23e-07} & \num{1.39e-06} & \num{7.94e-06} \\
    & \num{15} & \num{16384} & \num{6553600}
    & \num{2.05e-07} & \num{3.94e-07} & \num{2.26e-06} \\

    \bottomrule
  \end{tabular}
  \end{minipage}
\end{table}

\subsection{Exterior scattering and multiple right-hand sides}\label{subsec:ext_scatter}
The solver introduced in this manuscript is direct, so it can naturally handle multiple right hand sides. Table~\ref{tab:ext_scatter_results} reports the timings of simultaneously performing $500$ solves, corresponding to $500$ angles of incidence of plane waves. 

There are also many applications in which one is only interested in the scattered field at prescribed exterior locations. Examples include lens scattering~\cite{fuchs2008scattering,sanford1994scattering,mosallaei2001nonuniform}, microcavity emission~\cite{kim2016designing}, inverse
scattering~\cite{coltonkress2019,chew1990reconstruction,wei2019deep,ozmen2015comparing,borges2017high} and photonic crystal scattering~\cite{kristensen2010light,turduev2017photonic}. Significant effort has gone into designing specialized solvers to compute only the exterior scattered field (see, e.g., ~\cite{chew1992generalized,chew1993use,gillman2015spectrally}). In the present formulation, such a specialized solver is easily obtained by only retaining the root scatterer $\mtx{F}^1$. As Table~\ref{tab:ext_scatter_results} shows, this leads to large memory savings and significant speed-up.

Two types of exterior scattering are considered in Table~\ref{tab:ext_scatter_results}, which differ in the way the obtained outgoing field at $\partial\Omega$ is evaluated in the exterior of $\Omega$. After the outgoing trace is computed for each incoming field by applying $\mtx{F}^1$, both use the root level $\mtx{BIE}^{\rm(ext)}_{\rm solve}$ on $\partial\Omega$ followed by
\begin{enumerate}
    \item each exterior evaluation being performed separately, assuming the exterior evaluation points differ per incoming field (`direct'), or
    \item using a precomputed exterior evaluation matrix (EEM), mapping scattered fields at $\partial\Omega$ to the prescribed, shared, evaluation points (`from EEM').
\end{enumerate}
To be precise, assume the exterior points are given in advance. As in Section \ref{sec:disc_merge}, $\mtx{BIE}^{\rm(ext)}_{\rm eval}$ evaluates an exterior Hemholtz solution to these points. Then 
$\mtx{EEM} = \mtx{BIE}^{\rm(ext)}_{\rm eval} (\mtx{BIE}^{\rm(ext)}_{\rm solve} \,\mtx{F}^1)$, with $\mtx{BIE}^{\rm(ext)}_{\rm solve}$ the root-level exterior Brakhage–Werner solution operator.

In both cases, the exterior solve is orders of magnitude faster than the volume solve. Solving $500$ right-hand sides, even for the volume solve, is not significantly more expensive than solving a single right-hand side (see Table~\ref{tab:result_time}).
\begin{table}
  \centering
  \caption{Cavity scatterer ($\kappa = 201.06$, $p = 20$, $13$ levels): timings and
    storage for the volume solver and the exterior scatterer, for $500$ incidence
    directions. EEM denotes the exterior evaluation matrix; the error is the
    $\|\cdot\|_{\infty}$ difference from a fine solve.}
  \label{tab:ext_scatter_results}
  \begin{tabular}{lrrrrrrr}
    \toprule
      & \multicolumn{5}{c}{Time (s)} & & \\
    \cmidrule(lr){2-6}
      & Precomp. & Build & EEM & \multicolumn{2}{c}{Solve (500 rhs)} & Memory & Error \\
    \cmidrule(lr){5-6}
      &          &       &     & direct & from EEM & (GB) & \\
    \midrule
    Volume solver      & \num{46.7} & \num{54.9} & ---       & \num{14.92} & ---          & \num{10.6} & \num{1.81e-8} \\
    Exterior scatterer & \num{45.8} & \num{53.9} & \num{1.0} & \num{0.29}  & \num{2.1e-3} & \num{1.25} &         \num{1.72e-8}                       \\
    \bottomrule
  \end{tabular}
\end{table}

\begin{figure}
    \centering
    \includegraphics[width=0.8\linewidth]{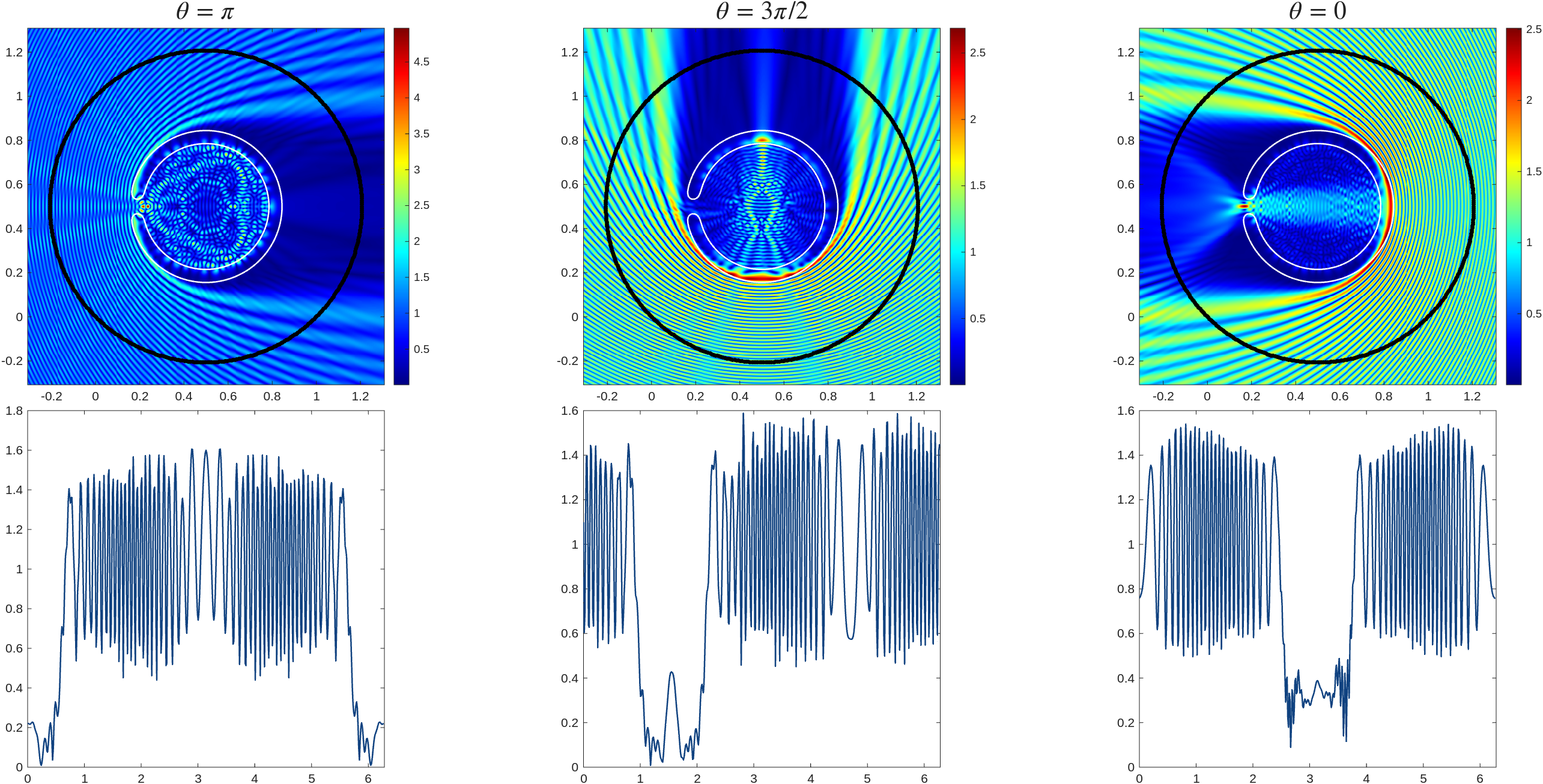}
    \caption{Magnitude of the total field $u_{\rm tot}$ for various incoming wave angles $\theta$. Top: global field. Bottom: Field measured at an exterior receiver array (indicated in black).}
    \label{fig:mult_scat_ext_fig}
\end{figure}

\section{Conclusions}
\label{sec:conc}

This paper describes a high-order accurate hierarchical direct solver for the two-dimensional Lippmann--Schwinger equation arising in acoustic scattering from variable media. The essential idea is to construct a reusable field-to-field scattering operator for the inhomogeneous medium, which would avoid a new global volume solve for every incident field. Each leaf subdomain is treated by a local Lippmann--Schwinger solve, while neighboring subdomains are merged recursively through boundary-only operations. As a result, volume solves are only operated on the finest level of the hierarchy, and the remainder of the computation is performed entirely in terms of boundary traces.

The resulting method naturally separates the computation into a reusable precomputation and an application stage. Once the scattering operator for a given medium, geometry, and frequency has been constructed, it can be applied to multiple incident fields without repeating the expensive volume discretization and hierarchical construction. Numerical experiments demonstrate high-order accuracy over a range of variable-media scattering problems, including smooth potentials as well as more challenging set-ups containing cavities, high-frequency problems, strong back-scattering, and other features that can make conventional iterative approaches difficult to use efficiently. The convergence results also illustrate the effect of the regularity of the scattering potential: smooth potentials exhibit rapid convergence, whereas potentials with limited regularity converge more gradually and may eventually be limited by the polynomial order used in the discretization.

Possible future work is to study the potential rank structure of the dense matrices arising throught the hierarchy. Although these matrices are dense, many of them have global low-rank structure, and the remaining matrices often exhibit hierarchical low-rank structure, such as hierarchically semiseparable (HSS) structure. In particular, exploiting low-rank and HSS representations could reduce both the storage requirements and the computational cost of the leaf solves, merging operations, and scattering-operator applications.

Overall, the numerical results demonstrate that the proposed hierarchical formulation provides an accurate and reusable direct solver for two-dimensional Lippmann--Schwinger scattering problems. Future work will focus on systematically exploiting the low-rank and HSS structure of the operators, with the goal of developing a compressed hierarchical implementation that retains the high-order accuracy while substantially reducing memory usage, construction cost, and application time.

\section*{Acknowledgements}

This work was supported by the Office of Naval Research (FP00011501), by the National Science Foundation (DMS-2608883 and DMS-2313434), and by the Simons Foundation. JB was supported in part by the Peter O'Donnell Jr.\ Postdoctoral Fellowship at the Oden Institute, and in part by the Swedish Research Council under grant no.\ 2021-06594 while in residence at Institut Mittag-Leffler in Djursholm, Sweden during the fall semester of 2025.

\clearpage
\bibliography{main_bib}

@article{2022_martinsson_lippmannschwinger,
  title={An accelerated, high-order accurate direct solver for the Lippmann--Schwinger equation for acoustic scattering in the plane},
  author={Gopal, Abinand and Martinsson, Per-Gunnar},
  journal={Advances in Computational Mathematics},
  volume={48},
  number={4},
  pages={1--31},
  year={2022},
  publisher={Springer}
}

@article{2019_ganesh_BIE_multiple_scattering,
title = "A fast high order algorithm for multiple scattering from large sound-hard three dimensional configurations",
journal = "Journal of Computational and Applied Mathematics",
volume = "362",
pages = "324 - 340",
year = "2019",
issn = "0377-0427",
doi = "https://doi.org/10.1016/j.cam.2018.10.053",
url = "http://www.sciencedirect.com/science/article/pii/S0377042718306617",
author = "M. Ganesh and S.C. Hawkins"
}

@article{corona2015n,
  title={An O (N) direct solver for integral equations on the plane},
  author={Corona, Eduardo and Martinsson, Per-Gunnar and Zorin, Denis},
  journal={Applied and Computational Harmonic Analysis},
  volume={38},
  number={2},
  pages={284--317},
  year={2015},
  publisher={Elsevier}
}

@book{2013_colton_kress,
  title={Integral equation methods in scattering theory},
  author={Colton, David and Kress, Rainer},
  volume={72},
  year={2013},
  publisher={SIAM}
}

@article{2013_martinsson_DtN_linearcomplexity,
author = {Gillman, A. and Martinsson, P.},
title = {A Direct Solver with $O(N)$ Complexity for Variable Coefficient Elliptic PDEs Discretized via a High-Order Composite Spectral Collocation Method},
journal = {SIAM Journal on Scientific Computing},
volume = {36},
number = {4},
pages = {A2023-A2046},
year = {2014},
doi = {10.1137/130918988},
URL = {http://dx.doi.org/10.1137/130918988},
eprint = {http://dx.doi.org/10.1137/130918988},
NOTE = {arXiv.org report \#1307.2665}
}

@article{2012_ganesh_multiple_scattering,
title = {A reduced basis method for electromagnetic scattering by multiple particles in three dimensions},
journal = {Journal of Computational Physics},
volume = {231},
number = {23},
pages = {7756-7779},
year = {2012},
issn = {0021-9991},
doi = {https://doi.org/10.1016/j.jcp.2012.07.008},
url = {https://www.sciencedirect.com/science/article/pii/S0021999112003671},
author = {M. Ganesh and J.S. Hesthaven and B. Stamm}
}

@article{2012_spectralcomposite,
title = {A direct solver for variable coefficient elliptic PDEs discretized via a composite spectral collocation method},
journal = {Journal of Computational Physics},
volume = {242},
number = {0},
pages = {460 - 479},
year = {2013},
note = {},
issn = {0021-9991},
doi = {10.1016/j.jcp.2013.02.019},
url = {http://www.sciencedirect.com/science/article/pii/S0021999113001320},
author = {P.G. Martinsson}
}

@article{2008_duan_rokhlin,
  author    = {Duan, Ran and
               Rokhlin, Vladimir},
  title     = {High-order quadratures for the solution of scattering problems
               in two dimensions},
  journal   = {J. Comput. Physics},
  volume    = {228},
  number    = {6},
  year      = {2009},
  pages     = {2152-2174},
  ee        = {http://dx.doi.org/10.1016/j.jcp.2008.11.033},
}

@book {2008_bebendorf_book,
    AUTHOR = {Bebendorf, Mario},
     TITLE = {Hierarchical matrices},
    SERIES = {Lecture Notes in Computational Science and Engineering},
    VOLUME = {63},
      NOTE = {A means to efficiently solve elliptic boundary value problems},
 PUBLISHER = {Springer-Verlag},
   ADDRESS = {Berlin},
      YEAR = {2008},
     PAGES = {xvi+290},
      ISBN = {978-3-540-77146-3},
   MRCLASS = {15-02 (35J40 65N22)},
  MRNUMBER = {2451321 (2009k:15001)},
MRREVIEWER = {Elena Pelican},
}

@article{2005_martinsson_fastdirect,
    AUTHOR = {Martinsson, P.G. and Rokhlin, V.},
     TITLE = {A fast direct solver for boundary integral equations in two dimensions},
   JOURNAL = {J. Comp. Phys.},
    VOLUME = {205},
    NUMBER = {1},
      YEAR = {2005},
     PAGES = {1--23},
}

@article {2003_hackbusch,
    AUTHOR = {Grasedyck, Lars and Hackbusch, Wolfgang},
     TITLE = {Construction and arithmetics of {$\mathcal{H}$}-matrices},
   JOURNAL = {Computing},
  FJOURNAL = {Computing. Archives for Scientific Computing},
    VOLUME = {70},
      YEAR = {2003},
    NUMBER = {4},
     PAGES = {295--334},
}

@article{2002_chen_direct_lippman_schwinger,
    AUTHOR = {Chen, Yu},
     TITLE = {A fast, direct algorithm for the {L}ippmann-{S}chwinger
              integral equation in two dimensions},
      NOTE = {Modeling and computation in optics and electromagnetics},
   JOURNAL = {Adv. Comput. Math.},
  FJOURNAL = {Advances in Computational Mathematics},
    VOLUME = {16},
      YEAR = {2002},
    NUMBER = {2-3},
     PAGES = {175--190},
      ISSN = {1019-7168},
   MRCLASS = {65R20 (78A45 78M25)},
  MRNUMBER = {1892248 (2003b:65139)},
       DOI = {10.1023/A:1014450116300},
       URL = {http://dx.doi.org/10.1023/A:1014450116300},
}

@article{rokhlin_1990,
AUTHOR = {Rokhlin, V.},
YEAR = {1990},
TITLE = {Rapid solution of integral equations of scattering
theory in two dimensions},
JOURNAL = {J. Comput. Phys.},
VOLUME = {86},
PAGES = {414--439}
}

@article{rokhlin1987,
    AUTHOR = {Greengard, L. and Rokhlin, V.},
     TITLE = {A fast algorithm for particle simulations},
   JOURNAL = {J. Comput. Phys.},
  FJOURNAL = {Journal of Computational Physics},
    VOLUME = {73},
      YEAR = {1987},
    NUMBER = {2},
     PAGES = {325--348},
      ISSN = {0021-9991},
     CODEN = {JCTPAH},
}

@article{2004_aguilar_chen_lippmann,
  author  = {Aguilar, J. C. and Chen, Yu},
  title   = {A high-order, fast algorithm for scattering calculation in two dimensions},
  journal = {Computers \& Mathematics with Applications},
  volume  = {47},
  number  = {1},
  pages   = {1--11},
  year    = {2004},
  doi     = {10.1016/S0898-1221(04)90001-6}
}

@article{2015_ying_lippmann_preconditioner,
  author  = {Ying, Lexing},
  title   = {Sparsifying Preconditioner for the {L}ippmann--{S}chwinger Equation},
  journal = {Multiscale Modeling \& Simulation},
  volume  = {13},
  number  = {2},
  pages   = {644--660},
  year    = {2015},
  doi     = {10.1137/140985147}
}

@article{2016_ambikasaran_lippmann,
  author  = {Ambikasaran, Sivaram and Borges, Carlos and Imbert-Gerard, Lise-Marie and Greengard, Leslie},
  title   = {Fast, Adaptive, High-Order Accurate Discretization of the {L}ippmann--{S}chwinger Equation in Two Dimensions},
  journal = {SIAM Journal on Scientific Computing},
  volume  = {38},
  number  = {3},
  pages   = {A1770--A1787},
  year    = {2016},
  doi     = {10.1137/15M102455X}
}

@article{2016_anand_pandey_lippmann,
  author  = {Anand, Akash and Pandey, Ambuj and Rathish Kumar, B. V. and Paul, Jagabandhu},
  title   = {An efficient high-order {N}ystr{\"o}m scheme for acoustic scattering by inhomogeneous penetrable media with discontinuous material interface},
  journal = {Journal of Computational Physics},
  volume  = {311},
  pages   = {258--274},
  year    = {2016},
  doi     = {10.1016/j.jcp.2016.01.028}
}

@article{2016_zepeda_zhao_lippmann,
  author  = {Zepeda-N{\'u}{\~n}ez, Leonardo and Zhao, Hongkai},
  title   = {Fast Alternating BiDirectional Preconditioner for the 2D High-Frequency {L}ippmann--{S}chwinger Equation},
  journal = {SIAM Journal on Scientific Computing},
  volume  = {38},
  number  = {5},
  pages   = {B866--B888},
  year    = {2016},
  doi     = {10.1137/16M1064660}
}

@article{2018_liu_ying_sparsify_sweep,
  author  = {Liu, Fei and Ying, Lexing},
  title   = {Sparsify and Sweep: An Efficient Preconditioner for the {L}ippmann--{S}chwinger Equation},
  journal = {SIAM Journal on Scientific Computing},
  volume  = {40},
  number  = {2},
  pages   = {B379--B404},
  year    = {2018},
  doi     = {10.1137/17M1132057}
}

@article{2018_chen_cho_cai_vie,
  author  = {Chen, Duan and Cho, Min Hyung and Cai, Wei},
  title   = {Accurate and Efficient {N}ystr{\"o}m Volume Integral Equation Method for Electromagnetic Scattering of 3-D Metamaterials in Layered Media},
  journal = {SIAM Journal on Scientific Computing},
  volume  = {40},
  number  = {1},
  pages   = {B259--B282},
  year    = {2018},
  doi     = {10.1137/16M110900X}
}

@article{2012_ganesh_hawkins_hiptmair_tmatrix,
  author  = {Ganesh, M. and Hawkins, S. C. and Hiptmair, R.},
  title   = {Convergence analysis with parameter estimates for a reduced basis acoustic scattering {T}-matrix method},
  journal = {IMA Journal of Numerical Analysis},
  volume  = {32},
  number  = {4},
  pages   = {1348--1374},
  year    = {2012},
  doi     = {10.1093/imanum/drr041}
}

@misc{chunkie,
  author = {Askham, Travis and Rachh, Manas and O'Neil, Michael and Hoskins, Jeremy and Fortunato, Dan and Jiang, Shidong},
  title = {chunkIE: a MATLAB integral equation toolbox},
  year = {2023},
  publisher = {GitHub},
  journal = {GitHub repository},
  howpublished = {\url{https://github.com/fastalgorithms/chunkie}}
}

@article{RCIP,
author = {Johan Helsing},
title = {{Solving Integral Equations on Piecewise Smooth Boundaries Using the RCIP Method: A Tutorial}},
volume = {2013},
journal = {Abstract and Applied Analysis},
number = {SI51},
publisher = {Hindawi},
pages = {1 -- 20},
year = {2013},
doi = {10.1155/2013/938167},
URL = {https://doi.org/10.1155/2013/938167}
}

@article{fuchs2008scattering,
  author  = {Fuchs, Benjamin and Palud, S{\'e}bastien and Le Coq, Laurent and
             Lafond, Olivier and Himdi, Mohamed and Rondineau, S{\'e}bastien},
  title   = {Scattering of Spherically and Hemispherically Stratified Lenses
             Fed by Any Real Source},
  journal = {IEEE Transactions on Antennas and Propagation},
  volume  = {56},
  number  = {2},
  pages   = {450--460},
  year    = {2008},
  doi     = {10.1109/TAP.2007.915458}
}

@article{sanford1994scattering,
  author  = {Sanford, John R.},
  title   = {Scattering by Spherically Stratified Microwave Lens Antennas},
  journal = {IEEE Transactions on Antennas and Propagation},
  volume  = {42},
  number  = {5},
  pages   = {690--698},
  year    = {1994}
}

@article{mosallaei2001nonuniform,
  author  = {Mosallaei, Hossein and Rahmat-Samii, Yahya},
  title   = {Nonuniform {Luneburg} and Two-Shell Lens Antennas: Radiation
             Characteristics and Design Optimization},
  journal = {IEEE Transactions on Antennas and Propagation},
  volume  = {49},
  number  = {1},
  pages   = {60--69},
  year    = {2001}
}

@article{kim2016designing,
  author  = {Kim, Yushin and Lee, Soo-Young and Ryu, Jung-Wan and Kim, Inbo and
             Han, Jae-Hyung and Tae, Heung-Sik and Choi, Muhan and Min, Bumki},
  title   = {Designing whispering gallery modes via transformation optics},
  journal = {Nature Photonics},
  volume  = {10},
  number  = {10},
  pages   = {647--652},
  year    = {2016},
  doi     = {10.1038/nphoton.2016.184}
}

@book{coltonkress2019,
  author    = {Colton, David and Kress, Rainer},
  title     = {Inverse Acoustic and Electromagnetic Scattering Theory},
  edition   = {4},
  series    = {Applied Mathematical Sciences},
  volume    = {93},
  publisher = {Springer},
  address   = {Cham},
  year      = {2019}
}

@article{chew1990reconstruction,
  author  = {Chew, Weng Cho and Wang, Y. M.},
  title   = {Reconstruction of Two-Dimensional Permittivity Distribution Using
             the Distorted {Born} Iterative Method},
  journal = {IEEE Transactions on Medical Imaging},
  volume  = {9},
  number  = {2},
  pages   = {218--225},
  year    = {1990}
}

@article{ozmen2015comparing,
  title={Comparing different ultrasound imaging methods for breast cancer detection},
  author={Ozmen, Neslihan and Dapp, Robin and Zapf, Michael and Gemmeke, Hartmut and Ruiter, Nicole V and Van Dongen, Koen WA},
  journal={IEEE transactions on ultrasonics, ferroelectrics, and frequency control},
  volume={62},
  number={4},
  pages={637--646},
  year={2015},
  publisher={IEEE}
}

@article{borges2017high,
  author  = {Borges, Carlos and Gillman, Adrianna and Greengard, Leslie},
  title   = {High Resolution Inverse Scattering in Two Dimensions Using
             Recursive Linearization},
  journal = {SIAM Journal on Imaging Sciences},
  volume  = {10},
  number  = {2},
  pages   = {641--664},
  year    = {2017},
  doi     = {10.1137/16M1093562}
}

@article{chew1993use,
  author  = {Chew, Weng Cho and Lu, Cai-Cheng},
  title   = {The use of {Huygens'} equivalence principle for solving the volume
             integral equation of scattering},
  journal = {IEEE Transactions on Antennas and Propagation},
  volume  = {41},
  number  = {7},
  pages   = {897--904},
  year    = {1993}
}

@article{chew1992generalized,
  title={A generalized recursive algorithm for wave-scattering solutions in two dimensions},
  author={Chew, Weng Cho and Gurel, Levent and Wang, Y-M and Otto, Gregory and Wagner, Robert L and Liu, Qing Huo},
  journal={IEEE transactions on microwave theory and techniques},
  volume={40},
  number={4},
  pages={716--723},
  year={1992},
  publisher={IEEE}
}

@article{gillman2015spectrally,
  author  = {Gillman, Adrianna and Barnett, Alex H. and Martinsson, Per-Gunnar},
  title   = {A spectrally accurate direct solution technique for frequency-domain
             scattering problems with variable media},
  journal = {BIT Numerical Mathematics},
  volume  = {55},
  number  = {1},
  pages   = {141--170},
  year    = {2015}
}

@article{kristensen2010light,
  title={Light propagation in finite-sized photonic crystals: multiple scattering using an electric field integral equation},
  author={Kristensen, Philip Tr{\o}st and Lodahl, Peter and M{\o}rk, Jesper},
  journal={Journal of the Optical Society of America B},
  volume={27},
  number={2},
  pages={228--237},
  year={2010},
  publisher={Optical Society of America}
}

@article{ying2004kernel,
  title={A kernel-independent adaptive fast multipole algorithm in two and three dimensions},
  author={Ying, Lexing and Biros, George and Zorin, Denis},
  journal={Journal of Computational Physics},
  volume={196},
  number={2},
  pages={591--626},
  year={2004},
  publisher={Elsevier}
}

@article{wei2019deep,
  author  = {Wei, Zhun and Chen, Xudong},
  title   = {Deep-Learning Schemes for Full-Wave Nonlinear Inverse Scattering
             Problems},
  journal = {IEEE Transactions on Geoscience and Remote Sensing},
  volume  = {57},
  number  = {4},
  pages   = {1849--1860},
  year    = {2019}
}

@article{turduev2017photonic,
  title={Photonic crystal based polarization insensitive flat lens},
  author={Turduev, Mirbek and Bor, E and Kurt, Hamza},
  journal={Journal of Physics D: Applied Physics},
  volume={50},
  number={27},
  pages={275105},
  year={2017},
  publisher={IOP Publishing}
}
\bibliographystyle{amsplain}

\clearpage

\begin{appendix}    

\section{Equivalent-source formulation}
\label{app:equivalent_source}

The field-to-field formulation in the main text represents both the incoming and outgoing
fields of a box by their Dirichlet traces. There is a closely related alternative in which
the incoming field is still represented by its trace $\phi_{\tau}$, but the outgoing field
is represented by an equivalent source density on the boundary. We record the formulation
here because its algebra is distinct, as it makes the propagation between boxes a
direct layer-potential evaluation, and it may prove useful in future variants of the method. The leaf solve is essentially unchanged; what differs is the representation of the
outgoing field and, consequently, the merge.

Throughout this appendix, $F^{\tau}$ denotes the equivalent-source scattering operator,
rather than the field-to-field scattering operator used in the main text.

\subsection{Local equivalent-source scattering operator}

Let $\Omega_{\tau}$ be a box and suppose that the incoming trace $\phi_{\tau}$ is given.
As in Section~\ref{sec:local_scattering}, first reconstruct the incoming field
$v_{\tau}=S_{\tau}\phi_{\tau}$ and solve the locally restricted Lippmann--Schwinger
equation. Equivalently, using (\ref{eq:def_source_operator}),
$$
q|_{\Omega_{\tau}}=Q^{\tau}\phi_{\tau}.
$$
The corresponding outgoing field is
$$
w_{\tau}(\pxx)=\int_{\Omega_{\tau}}G_{\kappa}(\pxx,\pyy)q(\pyy)\,d\pyy,
\qquad \pxx\in\mathbb{R}^{2}\setminus\overline{\Omega_{\tau}}.
$$
Instead of representing $w_{\tau}$ by its boundary trace, seek a density
$\mu_{\tau}$ on $\Gamma_{\tau}$ whose single-layer potential reproduces the same
exterior field:
\begin{equation}
\label{eq:eqsrc_equivalent_density}
\int_{\Omega_{\tau}}G_{\kappa}(\pxx,\pyy)q(\pyy)\,d\pyy
=
\int_{\Gamma_{\tau}}G_{\kappa}(\pxx,\pyy)\mu_{\tau}(\pyy)\,ds_{\pyy},
\qquad
\pxx\in\mathbb{R}^{2}\setminus\overline{\Omega_{\tau}}.
\end{equation}
Assuming that a unique equivalent density has been selected, define the
\textit{equivalent-source scattering operator}
\begin{equation}
\label{eq:eqsrc_local_scattering}
F^{\tau}:\phi_{\tau}\mapsto\mu_{\tau}.
\end{equation}
For the pure single-layer representation, uniqueness can fail at exceptional Helmholtz
frequencies. This issue can be avoided by restricting the box sizes, or by replacing the
single-layer representation by an appropriate combined-field representation. Since the
present appendix is intended only to record the alternative hierarchical formulation, we
do not pursue this question further.

\subsection{Merge of two adjacent boxes}

Consider the same geometry and boundary ordering as in Section~\ref{sec:mergetwo}:
$$
\partial\Omega_{\alpha}=\Gamma_{1}\cup\Gamma_{3},
\qquad
\partial\Omega_{\beta}=\Gamma_{2}\cup\Gamma_{4},
\qquad
\Gamma_{3}=\Gamma_{4},
$$
and
$$
\Gamma_{\tau}=\Gamma_{1}\cup\Gamma_{2}.
$$
The parent incoming trace $\phi_{\tau}$ is continued to the complete child boundaries
by the same universal operator $U_{\tau}$ used in the field-to-field formulation. Set
\begin{equation}
\label{eq:eqsrc_external_child_trace}
\varphi=U_{\tau}\phi_{\tau}
=
\left[\begin{array}{c}
\varphi_{1}\\
\varphi_{2}\\
\varphi_{3}\\
\varphi_{4}
\end{array}\right].
\end{equation}

Let
$$
\mu=
\left[\begin{array}{c}
\mu_{1}\\
\mu_{2}\\
\mu_{3}\\
\mu_{4}
\end{array}\right]
$$
denote the equivalent source densities on the two child boundaries, with
$$
\mu_{\alpha}=\left[\begin{array}{c}\mu_{1}\\\mu_{3}\end{array}\right],
\qquad
\mu_{\beta}=\left[\begin{array}{c}\mu_{2}\\\mu_{4}\end{array}\right].
$$
Split the child scattering operators according to their exterior and interface boundary
pieces,
$$
F^{\alpha}=\left[\begin{array}{cc}
F^{\alpha}_{11}&F^{\alpha}_{13}\\
F^{\alpha}_{31}&F^{\alpha}_{33}
\end{array}\right],
\qquad
F^{\beta}=\left[\begin{array}{cc}
F^{\beta}_{22}&F^{\beta}_{24}\\
F^{\beta}_{42}&F^{\beta}_{44}
\end{array}\right],
$$
and form the block operator
\begin{equation}
\label{eq:eqsrc_F}
F=
\left[\begin{array}{cccc}
F^{\alpha}_{11}&0&F^{\alpha}_{13}&0\\
0&F^{\beta}_{22}&0&F^{\beta}_{24}\\
F^{\alpha}_{31}&0&F^{\alpha}_{33}&0\\
0&F^{\beta}_{42}&0&F^{\beta}_{44}
\end{array}\right].
\end{equation}
Thus $F$ maps the total incoming traces on the two child boundaries to the corresponding
equivalent source densities.

The field generated by the sources on one child provides part of the incoming field on
its sibling. For example,
\begin{equation}
\label{eq:eqsrc_G_pairwise}
\left[\begin{array}{c}
\xi_{1}\\
\xi_{3}
\end{array}\right]
=
\left[\begin{array}{cc}
G_{12}&G_{14}\\
G_{32}&G_{34}
\end{array}\right]
\left[\begin{array}{c}
\mu_{2}\\
\mu_{4}
\end{array}\right],
\qquad
\left[\begin{array}{c}
\xi_{2}\\
\xi_{4}
\end{array}\right]
=
\left[\begin{array}{cc}
G_{21}&G_{23}\\
G_{41}&G_{43}
\end{array}\right]
\left[\begin{array}{c}
\mu_{1}\\
\mu_{3}
\end{array}\right].
\end{equation}
Here $\xi$ denotes the sibling-generated incoming trace. Each operator $G_{ij}$ is simply
a layer-potential evaluation operator: it maps a single-layer density on $\Gamma_j$ to
the value of the resulting field on $\Gamma_i$. For instance,
$$
[G_{12}\mu_{2}](\pxx)=\int_{\Gamma_{2}}G_{\kappa}(\pxx,\pyy)\mu_{2}(\pyy)\,ds_{\pyy},
\qquad \pxx\in\Gamma_{1}.
$$
Combining the two relations in (\ref{eq:eqsrc_G_pairwise}) gives
\begin{equation}
\label{eq:eqsrc_G}
\xi=G_{\tau}\mu,
\qquad
G_{\tau}=
\left[\begin{array}{cccc}
0&G_{12}&0&G_{14}\\
G_{21}&0&G_{23}&0\\
0&G_{32}&0&G_{34}\\
G_{41}&0&G_{43}&0
\end{array}\right].
\end{equation}
The merge equilibrium equation is therefore
\begin{equation}
\label{eq:eqsrc_equilibrium}
\mu=F\bigl(\varphi+G_{\tau}\mu\bigr).
\end{equation}
Using (\ref{eq:eqsrc_external_child_trace}), we obtain
\begin{equation}
\label{eq:eqsrc_child_sources}
\mu=(I-FG_{\tau})^{-1}F U_{\tau}\phi_{\tau}.
\end{equation}

It remains to replace the four child densities in $\mu$ by one equivalent density on the
parent boundary. Write
$$
\mu_{\tau}=
\left[\begin{array}{c}
\widehat\mu_{1}\\
\widehat\mu_{2}
\end{array}\right],
$$
where $\widehat\mu_{1}$ and $\widehat\mu_{2}$ are supported on $\Gamma_{1}$ and
$\Gamma_{2}$, respectively. Define the equivalent-source assembly operator $V_{\tau}$ by
\begin{equation}
\label{eq:eqsrc_V_block}
\left[\begin{array}{c}
\widehat\mu_{1}\\
\widehat\mu_{2}
\end{array}\right]
=
\underbrace{\left[\begin{array}{cccc}
I&0&V_{13}&V_{14}\\
0&I&V_{23}&V_{24}
\end{array}\right]}_{=:V_{\tau}}
\left[\begin{array}{c}
\mu_{1}\\
\mu_{2}\\
\mu_{3}\\
\mu_{4}
\end{array}\right].
\end{equation}
The blocks $V_{ij}$ are chosen so that the density $\mu_{\tau}=V_{\tau}\mu$ generates
the same field outside $\Omega_{\tau}$ as the two child single-layer potentials. Equivalently,
\begin{multline}
\label{eq:eqsrc_parent_representation}
\int_{\Gamma_{\tau}}G_{\kappa}(\pxx,\pyy)\mu_{\tau}(\pyy)\,ds_{\pyy}
=\\
\int_{\partial\Omega_{\alpha}}G_{\kappa}(\pxx,\pyy)\mu_{\alpha}(\pyy)\,ds_{\pyy}
+
\int_{\partial\Omega_{\beta}}G_{\kappa}(\pxx,\pyy)\mu_{\beta}(\pyy)\,ds_{\pyy},
\qquad
\pxx\in\mathbb{R}^{2}\setminus\overline{\Omega_{\tau}}.
\end{multline}
Combining (\ref{eq:eqsrc_child_sources}) and (\ref{eq:eqsrc_V_block}) gives the parent
scattering operator
\begin{equation}
\label{eq:eqsrc_merge}
\boxed{
F^{\tau}=V_{\tau}(I-FG_{\tau})^{-1}F U_{\tau}.
}
\end{equation}
Thus the equivalent-source formulation has the same basic hierarchical structure as the
field-to-field method: solve a local multiple-scattering equilibrium between the children,
and then replace the pair by a single scattering operator for the parent.

\subsection{Structure of the merge}

The product $FG_{\tau}$ that appears in (\ref{eq:eqsrc_child_sources}) has a particularly
useful block structure. From (\ref{eq:eqsrc_F}) and (\ref{eq:eqsrc_G}),
\begin{equation}
\label{eq:eqsrc_FG}
FG_{\tau}=
\left[\begin{array}{cccc}
0 & F^{\alpha}_{11}G_{12}+F^{\alpha}_{13}G_{32} & 0 & F^{\alpha}_{11}G_{14}+F^{\alpha}_{13}G_{34}\\
F^{\beta}_{22}G_{21}+F^{\beta}_{24}G_{41} & 0 & F^{\beta}_{22}G_{23}+F^{\beta}_{24}G_{43} & 0\\
0 & F^{\alpha}_{31}G_{12}+F^{\alpha}_{33}G_{32} & 0 & F^{\alpha}_{31}G_{14}+F^{\alpha}_{33}G_{34}\\
F^{\beta}_{42}G_{21}+F^{\beta}_{44}G_{41} & 0 & F^{\beta}_{42}G_{23}+F^{\beta}_{44}G_{43} & 0
\end{array}\right].
\end{equation}
This identity is useful because almost every term in (\ref{eq:eqsrc_FG}) contains at least
one factor that is numerically low rank. The principal exceptions are the terms
$F^{\alpha}_{33}G_{34}$ and $F^{\beta}_{44}G_{43}$, which couple the two copies of the
common interface. These blocks are not globally low rank, but they have multilevel
rank structure. Thus, just as in the field-to-field formulation, the merge contains
considerable exploitable structure even though the present appendix does not attempt to
use it.

\begin{remark}[Comparison with the field-to-field formulation]
The two formulations use the same incoming trace and the same local
Lippmann--Schwinger solves. Their difference is entirely in the representation of the
outgoing field. In the formulation used in the main text, the outgoing state is the
boundary trace $\psi_{\tau}$. This avoids equivalent-source reconstruction and makes the
coupling across the common interface particularly simple: two blocks of $E_{\tau}$ are
identities. In the present formulation, the outgoing state is the single-layer density
$\mu_{\tau}$. Propagation between children is then given directly by the universal
single-layer evaluation operator $G_{\tau}$, but the near-field blocks $G_{34}$ and
$G_{43}$ require singular or nearly singular layer-potential evaluation, and the parent
merge requires the additional equivalent-source reconstruction $V_{\tau}$. These
considerations are the main reasons that the field-to-field formulation was preferred in
the implementation described in this paper.
\end{remark}

\begin{remark}[Universal operators]
The interaction operator $G_{\tau}$ and the assembly operator $V_{\tau}$ are independent
of the scattering potential $b$, as is the incoming continuation operator $U_{\tau}$. They
depend only on the geometry and on the free-space Helmholtz operator, and can therefore
be precomputed. The separated-boundary blocks of $G_{\tau}$ are smoothing and highly
compressible. The interface blocks $G_{34}$ and $G_{43}$ require near-field evaluation and,
at high accuracy, a multilevel representation. These observations provide a possible
alternative route to accelerating the method if the equivalent-source formulation is
revisited in future work.
\end{remark}

\end{appendix}

\end{document}